\documentclass[a4paper,11pt,twoside]{article}

\usepackage{latexsym,amsmath,amssymb,amsfonts,amsbsy,fancyhdr}
\usepackage{theorem}
\usepackage{amscd}
\usepackage{enumerate}
\usepackage{mathrsfs}
\usepackage{dsfont}
\usepackage[active]{srcltx}
\usepackage[colorlinks=true, pdfstartview=FitV, linkcolor=black, citecolor=black, urlcolor=black]{hyperref}
\usepackage[dvips]{graphicx}
\usepackage{epsf}
\usepackage{nicefrac}
\usepackage{marvosym}
\usepackage{subfigure}

\newcommand{\R}{{\mathbb R}}
\newcommand{\rp}{\R^+}
\newcommand{\rpn}{\R^+_0}

\newcommand{\N}{{\mathbb N}}

\newcommand{\dd}{\text{d}}

\newcommand{\comment}[1]{\marginpar{\textcolor{red}{#1}}}

\newtheorem{theorem}{Theorem}

\newtheorem{proposition}[theorem]{Proposition}
\newtheorem{lemma}[theorem]{Lemma}
\theorembodyfont{\rm}

\newtheorem{example}[theorem]{Example}

\newtheorem{remark}[theorem]{Remark}
\newenvironment{rem}{\begin{remark}}{\hfill$\lozenge$\end{remark}}
\newtheorem{proof}{Proof}

\newenvironment{pf}{\begin{proof}}{\hfill$\square$\end{proof}}
\newtheorem{hypothesis}{Hypothesis}

\allowdisplaybreaks

\DeclareMathAlphabet{\mathcal}{OMS}{cmsy}{m}{n}

\begin{document}
\title{Asymptotic properties of the process counted with a random characteristic in the context of fragmentation processes}
\author{Robert Knobloch\thanks{Institut f\"ur Mathematik, FB 12, Goethe-Universit\"at Frankfurt am Main, 60054 Frankfurt am Main, Germany
\newline
e-mail: knobloch@math.uni-frankfurt.de}}
\maketitle

\makeatletter
   \let\my@makefnmark\@makefnmark
   \renewcommand{\@makefnmark}{}
   \makeatother
   
\begin{abstract}
In this paper we prove a strong law of large numbers and its $\mathscr L^1$-convergence counterpart for the process counted with a random characteristic in the context of self-similar fragmentation processes. This result extends a somewhat analogical result by Nerman for general branching processes to fragmentation processes. In addition, we apply the general result of this paper to a specific example that in particular extends a limit theorem, concerning the fragmentation energy, by Bertoin and Mart\'inez  from $\mathscr L^1$-convergence to almost sure convergence. Our approach treats fragmentation processes with an infinite dislocation measure directly, without using a discretisation method. Moreover, we obtain a result regarding the asymptotic behaviour of the empirical mean associated with some stopped fragmentation process.
\end{abstract}




\section{Introduction}
In the present paper we introduce random characteristics associated with self-similar fragmentation processes. In the setting of general branching processes (also known as Crump-Mode-Jagers processes) random characteristics were considered for example in \cite{91}, \cite{80}, \cite{88} as well as \cite{82} and  recently by Vatutin \cite{Vat11}. The motivation that led to this work is multifaceted and stems  from the literature on general branching processes and fragmentation processes respectively.  Let us dwell for a little while on outlining some related results in the literature which motivated our considerations. In \cite{80} Nerman proved various limit results for general branching processes counted with a random characteristic. There he managed to prove $\mathscr L^1$-convergence (cf. \cite[Corollary~3.3]{80}) as well as almost sure convergence (see \cite[Theorem~5.4]{80} and \cite[Theorem~6.3]{80}).  By setting {\it birth time of a particle} to be the negative logarithm of the {\it size of a block} one can  interpret, for the topic of this paper,    fragmentations with a finite dislocation measure as special cases of general branching processes. Indeed, this is possible, since the time-parameter of the fragmentation process does not affect the results.  Under suitable assumptions, which ensure that the conditions of Nerman's results are satisfied, one does then obtain respective results for fragmentation processes. This observation lies at the heart of \cite[Theorem~1, Corollary~1]{84}. However, with regard to fragmentation processes an interesting question is whether such results can be obtained also in the general case if the dislocation measure is  $\sigma$-finite rather than finite. In \cite{84} Bertoin and Mart\'inez proved $\mathscr L^1$-convergence (cf. \cite[Theorem~2]{84}) in the dissipative case, as well as $\mathscr L^2$-convergence (see \cite[Theorem~3]{84}) in the conservative case, for some functional (cf. Section~\ref{ss.erc} of the present paper) that can be considered as a special case of a process counted with a random characteristic. Their method of proving their Theorem~2 is based on using Theorem~7.3 of  \cite{88}, which extends Corollary~3.3 of \cite{80} to random characteristics endowed with a {\it life career}, and the aforementioned interpretation of fragmentations in the case of a finite dislocation measure as well as a discretisation method.  The proof of Theorem~3 in \cite{84} is based on moment calculations rather than using a discretisation technique. However, neither of the two methods can be reasonably modified to yield almost sure convergence. In Theorem~1 of \cite{HKK10} almost sure convergence of the process counted with the random characteristic considered in Corollaries~1 and 2 of \cite{84}  is proven. The approach in  \cite{HKK10} treats general dissipative fragmentation processes directly, without resorting to known results on general branching processes or fragmentations with a finite dislocation measure.  This raises the question whether by means of such a direct approach the $\mathscr L^1$-convergence of Theorem~2 in \cite{84} can be proven to  hold also almost surely and whether similar results hold for a reasonably large class of random characteristics in the spirit of \cite{80}. The present paper answers both of these questions affirmatively.

The outline of this paper is as follows.  In Section~\ref{s.hfp} we provide a  concise introduction to the theory of fragmentation processes and the associated process stopped at a stopping line. Subsequently, in Section~\ref{s.mr} we introduce random characteristics, as well as the process counted with them, and present our main result. In Section~\ref{ss.erc} we consider two applications of our strong law of large numbers. These examples are related to the discussion  above and provide some motivation for the general result. In Section~\ref{ss.p.{t.2.1}} we present the proof of our main result and Section~\ref{s.app} is devoted to proving a crucial ingredient of that proof. 

Throughout this paper we consider a probability space $(\Omega,\mathscr F,\mathbb P)$ on which the fragmentation process as well as all the other random objects are defined.

\section{Fragmentation processes}\label{s.hfp}
In this section we provide a brief introduction to fragmentation processes and associated  processes that are stopped at a stopping line. Moreover, we introduce two families of additive martingales. 

Let $\mathcal{P}$ be the space of partitions $\pi=(\pi_n)_{n\in\N}$ of $\N$, where the blocks of $\pi$ are ordered by their least element such that $\inf(\pi_i)<\inf(\pi_j)$ if $i<j$, where $\inf(\emptyset):=\infty$. This paper is concerned with a $\mathcal P$-valued fragmentation process $\Pi:= (\Pi(t))_{t\in\rpn}$, where $\Pi(t) = (\Pi_n(t))_{n\in\N}$. $\mathcal{P}$-valued fragmentations are exchangeable Markov processes that were introduced in  \cite{85} in the homogenous case and were extended to the self-similar setting in \cite{105}. For a comprehensive treatise on fragmentation processes we refer to the monograph  \cite{92}. Let $(\mathscr F_t)_{t\in\rpn}$ be the filtration generated by $\Pi$. 

It is known from \cite{105} that the distribution of $\Pi$ is determined by some $\alpha\in\R$ (the index of self-similarity; $\alpha=0$ corresponding to the homogenous case), a constant $c\in\rpn$ (the rate of erosion) and a measure $\nu$ (the so-called dislocation measure that determines the jumps of $\Pi$) on the infinite simplex
 \[
\mathcal S:=\left\{{\bf s}:=(s_n)_{n\in\N}:s_1\ge s_2\ldots\ge0,\,\sum_{n\in\N}s_n\le1\right\},
 \] 
such that  $\nu(\{(1,0,\ldots)\})=0$ as well as 
\[
\int_{\mathcal{S}}(1-s_1)\nu(d{\bf s})<\infty.
\]
Below we shall need the following constant 
\[
\underline{ p}: = \inf\left\{  p\in \mathbb{R} : \int_{\mathcal S} \left| 1- \sum_{n\in\N} s_n^ {1+p} \right| \nu(d{\bf s}) <\infty\right\}\in(-1, 0]
\]
as well as the increasing and concave function $\Phi:(\underline p,\infty)\to\R$, given by
\[
\Phi(p)=\int_{\mathcal S}\left(1-\sum_{n\in\N}s_n^{1+p}\right)\nu(\dd{\bf s})
\]
for every $(\underline p,\infty)$. The function $\Phi$ plays a crucial role in the theory of fragmentations, since it turns out to be the Laplace exponent of a killed subordinator defined via a fragmentation process. 

{\it Throughout this paper we consider a self-similar fragmentation process $\Pi$  that satisfies
\begin{equation}\label{e.cE}
c=0,\qquad\nu(\mathcal S)=\infty\qquad\text{and}\qquad\nu({\bf s}\in\mathcal S:s_2=0)<\infty
\end{equation}} 
as well as the following hypothesis which is often referred to as {\it Malthusian hypothesis}.
\begin{hypothesis}\label{h.1}
There exists a $p^*\in(\underline p,0]$ such that $\Phi(p^*)=0$.
\end{hypothesis}
The constant $p^*$ is often called {\it Malthusian parameter} 
and if $\nu$ is finite this definition  coincides with the  definition of the Malthusian parameter in the context of general branching processes. 

The process $\Pi$ is said to be {\it conservative} if $\nu(\sum_{n\in\N}s_n<1)=0$, i.e. if there is no loss of mass by sudden dislocations, and {\it dissipative} otherwise. In this paper we allow for both of these cases. 

Most proofs of this paper shall be established under the assumption that  $\Pi$ is homogenous, i.e.  $\alpha=0$, since many important concepts in the theory of fragmentation processes only hold in the homogenous setting. That is to say, we will first prove our main result for the homogenous case and then extend it to the  self-similar setting with $\alpha\in\R$. The crucial point in this regard is that every self-similar fragmentation process is a time-changed homogenous fragmentation process (see  Theorem~3~(i) of \cite{105}). In this spirit, for the remainder of this section we assume that $\alpha=0$.

We shall need the  exchangeable partition measure $\mu$ on $\mathcal{P}$ given by
\[
\mu(d\pi) = \int_{\mathcal{S}}\varrho_{\bf s}(d\pi)\nu(d{\bf s}),
\]
where $\varrho_{\bf s}$ is the law of Kingman's paint-box based on ${\bf s}\in\mathcal S$. In \cite{85} Bertoin showed that the homogenous  fragmentation process $\Pi$ is characterised by a Poisson point process. More precisely, there exists a $\mathcal P\times\mathbb N$-valued Poisson point process $(\pi(t),\kappa(t))_{t\in\rpn}$\label{p.pi_t} with characteristic measure $\mu\otimes\sharp$, where $\sharp$ denotes the counting measure on $\N$,  such that $\Pi$ changes state precisely at the times $t\in\rpn$ for which an atom $(\pi(t),\kappa(t))$ occurs in $(\mathcal P\setminus(\N,\emptyset,\ldots))\times\mathbb N$. At such a time $t\in\rpn$ the sequence $\Pi(t)$ is obtained from $\Pi(t-)$ by replacing its $\kappa(t)$-th term, $\Pi_{\kappa(t)}(t-)\subseteq\N$, with the restricted partition $\pi(t)|_{\Pi_{\kappa(t)}(t-)}$ and reordering the terms such that the resulting partition of $\N$ is an element of $\mathcal P$. We denote the random jump times of $\Pi$, i.e. the times at which the abovementioned Poisson point process has an atom in $(\mathcal P\setminus(\N,\emptyset,\ldots))\times\mathbb N$, by $(t_i)_{i\in\mathcal I}$, where the index set $\mathcal I\subseteq\rpn$ is  countably infinite.  

Moreover, by exchangeability, the limits
\[
|\Pi_n(t)|:=\lim_{k\to\infty}\frac{\sharp(\Pi_n(t)\cap\{1,\ldots,k\})}{k},
\]
referred to as {\it asymptotic frequencies}, exist $\mathbb P$-a.s. simultaneously for all $t\in\rpn$ and $n\in\N$. Let us point out that the concept of asymptotic frequencies provides us with a notion of {\it size} for the blocks of a $\mathcal P$-valued fragmentation process. In \cite{85} Bertoin showed that the process $(-\ln(|\Pi_1(t)|))_{t\in\rpn}$ is a killed subordinator with Laplace exponent $\Phi$, a fact we shall make use of below. 

\begin{figure}[htb]
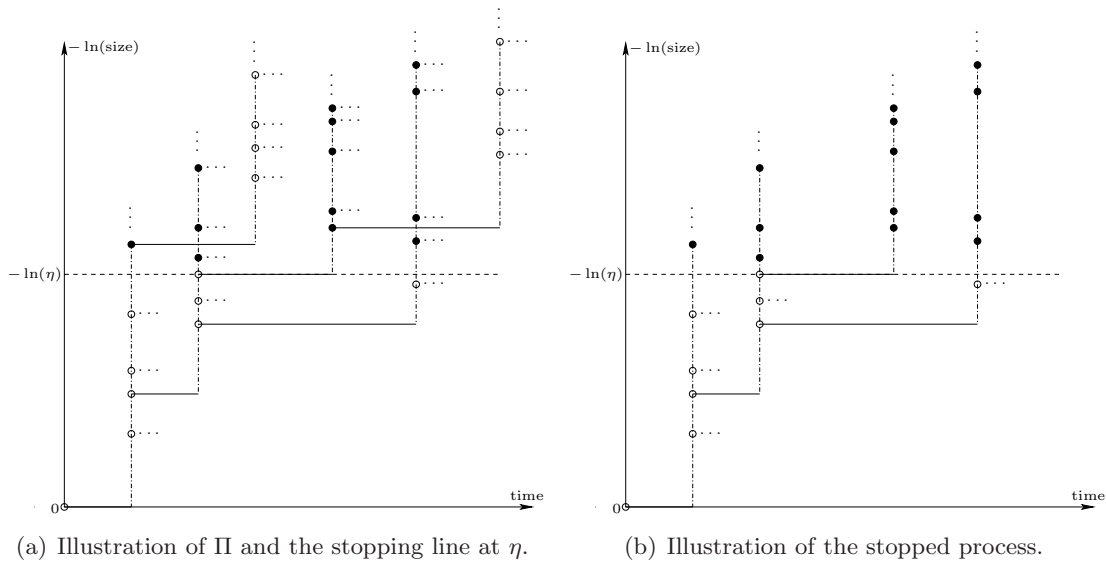
\centering
\subfigure[Illustration of $\Pi$ and the stopping line at $\eta$.\label{f.1.2a}]{\resizebox{7cm}{!}{\input{frag_stopping_line_50_eta.pstex_t}}}\quad
\subfigure[Illustration of the  stopped process.\label{f.1.2b}]{\resizebox{7cm}{!}{\input{frag_stopping_line_50b_eta.pstex_t}}}\quad
\caption[Stopped fragmentation process]{Illustration (a) depicts a realisation of a  fragmentation process $\Pi$ (with finite dislocation measure) and the stopping line  given by the first passage of the block sizes below some $\eta\in(0,1]$. In (b)  the fragmentation process  which is stopped at this stopping line is illustrated. The black dots indicate the blocks in $\lambda_\eta$, since their sizes are smaller than $\eta$ and they result from the dislocation of blocks with size greater than or equal to $\eta$. 
}\label{f.1.2}
\end{figure}

The main  property of fragmentation processes 
 is the {\it (strong) fragmentation property}, which is the analogue of the branching property of branching processes. Roughly speaking, this property says that given a configuration of the process at some (stopping) time, the  further evolution of each block is governed by an independent copy of the original process. Moreover, the same holds true if we replace the stopping time by a stopping line, cf. Definition~3.4 and Lemma~3.14 both in \cite{92}, in which case we refer to it as the {\it extended fragmentation property}. Speaking of stopping lines, let us now introduce stopped fragmentations which are obtained from a fragmentation process by stopping the evolution of a block once it has reached a given stopping line. Here we are interested in the stopping line that corresponds to the first blocks, in their respective ``line of descent'', of size less than some given $\eta\in(0,1]$. The  process $(\lambda_\eta)_{\eta\in(0,1]}$ consisting for each $\eta\in(0,1]$ of the sizes of the blocks at the terminal state of the fragmentation  stopped at the stopping line associated with $\eta$   is then given by $\lambda_\eta:=(\lambda_{\eta,k})_{k\in\N}$, where $\lambda_{\eta,k}$ refers  to the asymptotic frequency of the $k$-th largest block at the terminal state of the stopped process. In addition, we denote by $(\mathscr H_\eta)_{\eta\in(0,1]}$ the filtration generated by $(\lambda_\eta)_{\eta\in(0,1]}$. For an illustration of 
these concepts, see Figure~\ref{f.1.2}.
\\

In what follows we shall make use of two additive unit-mean martingales. Let  $\bar p>p^*$ be  given by Lemma~1 in \cite{89}. The first martingale, $M(p):=(M_t(p))_{t\in\rpn}$, is given by
\[
M_t(p):=\sum_{k\in\N}|\Pi_{k}(t)|^{1+p}e^{\Phi(p)t}
\]
for every $p\in(\underline p,\bar p)$ and $t\in\rpn$. Let us point out that $M(p)$ is the analogue of Biggins'  additive martingale for branching random walks (cf. \cite{Big77}). These martingales  appear frequently in the literature on branching processes and fragmentation processes (see e.g. \cite{84} where it was used in a context related to the present paper). In the spirit of \cite{80}  and \cite{HKK10} we shall need a second  martingale which is defined in terms of $(\lambda_\eta)_{\eta\in(0,1]}$ and turns out to be related to $M(p^*)$. More precisely, consider the process $(\Lambda_\eta(p^*))_{\eta\in(0,1]}$   defined by
\[
\Lambda_\eta(p^*):=\sum_{k\in\N}\lambda_{\eta,k}^{1+p}
\]
for each $\eta\in(0,1]$
The martingale property of $(\Lambda_\eta(p^*))_{\eta\in(0,1]}$ was established in Lemma~1 of \cite{HKK10}. In particular, there it was shown that 
\[
\Lambda_\eta(p^*)=\mathbb E\left(\left.M_\infty(p^*)\right|\mathscr H_\eta\right)
\]
holds for all $\eta\in(0,1]$, where $M_\infty(p^*):=\lim_{t\to\infty}M_t(p^*)$. As we will see in the next section, the limit that  appears in our main result turns out to be a deterministic constant times the almost sure martingale limit 
\begin{equation}\label{e.Lambda_0}
\Lambda_0(p^*):=\lim_{\eta\downarrow0}\Lambda_\eta(p^*).
\end{equation}

\section{Main result}\label{s.mr}


Recall that the present paper is concerned with a self-similar fragmentation process $\Pi$ that satisfies (\ref{e.cE}) as well as Hypothesis~\ref{h.1}. 

Throughout the present paper let $(\phi(x,\pi))_{x\in\rpn}$, $\pi\in\mathcal P$, be an $\R$-valued  stochastic processes,  with $\phi(x,\pi)=\phi(0,\pi)=0$ $\mathbb P$-a.s. for all $x>1$, which 
has c\`adl\`ag paths $\mathbb P$-a.s. and is independent of $\Pi$. In addition, assume that $\phi(\cdot,(1,0,\ldots))=0$. In the context of this paper we refer to the random function $\phi:\rpn\times\mathcal P\times\Omega\to\R$ as {\it random characteristic}.

Recall that if $\Pi$ is homogenous (i.e. $\alpha=0$), then it is determined by the Poisson point process $(\pi(t),\kappa(t))_{t\in\rpn}$. If $\alpha\ne0$, then we cannot resort to such a Poissonian structure of $\Pi$. However, in the light of Theorem~3~(i) in \cite{105} we still have that every jump of a block of $\Pi$ is determined by an independent copy of $\pi(1)$. Hence, we denote by $(\pi(t,k))_{t\in\rpn,k\in\N}$ the $\mathcal P$-valued random field that determines the jumps of $\Pi$ such that at each jump time $t\in\rpn$ of some block $\Pi_k(t-)$ this block fragments into the sequence $(\pi_n(t,k)|_{\Pi_k(t-)})_{n\in\N}$. If the block $\Pi_k(t-)$ does not jump at time $t$, then $\pi(t,k)=(1,0,\ldots)$.

In this work we are interested in the process $(Z^\phi_\eta)_{\eta\in(0,1]}$, often referred to as the {\it process counted with the characteristic $\phi$}, 
given by
\[
Z^\phi_\eta:=\sum_{(t,k)\in\mathcal T\times\N}\phi^{(t,k)}\left(\frac{\eta}{|\Pi_k(t-)|},\pi(t,k)\right)\mathds1_{\{|\Pi_k(t-)|\ne0\}}
\]
for every $\eta\in(0,1]$, where the $\phi^{(t,k)}$ are independent copies of $\phi$ and where
\[
\mathcal T:=\{t\in\rpn:\Pi(t)\ne\Pi(t-)\}
\] 
denotes the countably infinite set of jump times of $\Pi$.  In the light of the Poissonian structure of $\Pi$ in the homogenous setting we have
\[
Z^\phi_\eta=\sum_{i\in\mathcal I}\phi^{(i)}\left(\frac{\eta}{|\Pi_{\kappa(t_i)}(t_i-)|},\pi(t_i)\right)\mathds1_{\left\{|\Pi_{\kappa(t_i)}(t_i-)|\ne0\right\}}
\]
for each $\eta\in(0,1]$ if $\Pi$ is homogenous, where the $\phi^{(i)}$ are independent copies of $\phi$. 
In order to avoid the indicator function $\mathds1_{\{|\Pi_{\kappa(t_i)}(t_i-)|\ne0\}}$, we adopt the convention
\[
\phi\left(\frac{a}{0},\pi\right):=0
\]
for all $a\in\rpn$ and $\pi\in\mathcal P$.

The main goal of this paper is to prove  an extension of \cite[Theorem~5.4]{80} to fragmentation processes that in particular yields an extension of \cite[Theorem~1]{84} to almost sure convergence. As we shall show in Section~\ref{ss.erc}, such an extension of \cite[Theorem~1]{84} follows from the general result the present paper is concerned with. Generally speaking, we aim at deriving some asymptotic properties of $Z^\phi_\eta$ as $\eta\downarrow0$. More precisely, we prove that asymptotically, as $\eta\downarrow0$, $Z^\phi_\eta$ behaves like the product of its expectation and  $\Lambda_\eta(p^*)$. Indeed, the main result of this paper is strong law of large numbers. Denote by $\Phi'$ the derivative of $\Phi$ and recall $\lambda_0(p^*)$ given by (\ref{e.Lambda_0}).
\begin{theorem}\label{t.2.1}
Let $\phi$ be a random characteristic such that 
\begin{equation}\label{e.assumption_1}
\int_{\mathcal P}\mathbb E\left(\sup_{\eta\in(0,1]}\eta^{(1+p^*+\beta)}\phi\left(\eta,\pi\right)\right)\mu(\dd\pi)<\infty\quad\text{and}\quad\limsup_{\eta\downarrow0}\int_{\mathcal P}\eta^{1+\tilde p}\mathbb E\left(\phi\left(\eta,\pi\right)\right)\mu(\dd\pi)<\infty
\end{equation}
hold for  all $\beta>0$ and some $\tilde p\in(\underline p,p^*)$.
Then 
\[
\eta^{1+p^*}Z^\phi_\eta\to \frac{\Lambda_0(p^*)}{\Phi'(p^*)}\int_{(0,1)}\mathbb E\left(\sum_{k\in\N}|\Pi_{k}(t)|^{1+p^*}\int_\mathcal P\int_{(0,1]}\rho^{p^*}\phi\left(\rho,\pi\right)\dd\rho\mu(\dd\pi)\right)\dd t
\]
$\mathbb P$-a.s. and in $\mathscr L^1(\mathbb P)$ as $\eta\downarrow0$.
\end{theorem}
We give the proof of Theorem~\ref{t.2.1} in Section~\ref{ss.p.{t.2.1}}. 
Of course, the statement of Theorem~\ref{t.2.1}  remains true if we replace (\ref{e.assumption_1}) by the stronger assumption
\begin{equation}\label{e.assumption_1_strong}
\int_\mathcal P\mathbb E\left(\sup_{\eta\in(0,1]}\eta^{1+\tilde p}\phi\left(\eta,\pi\right)\right)\mu(\dd\pi)<\infty
\end{equation}
for some $\tilde p\in(\underline p,p^*)$. 
Recall that our goal was to extend \cite[Theorem~5.4]{80} to fragmentation processes with an infinite dislocation measure. In this spirit, notice that  (\ref{e.assumption_1_strong}) is in some sense a natural analogue for fragmentations  of Condition~5.2 in \cite{80} that appears in \cite[Theorem~5.4]{80}. Admittedly, (\ref{e.assumption_1_strong}) is actually more restrictive than Nerman's Condition~5.2 as it requires the integrability with respect to $\mu$, but in contrast to general branching processes such a requirement is in general necessary for fragmentation processes. 

A simple example of a random characteristic is 
\[
\phi(x,\pi):=\mathds1_{\{x\in(a,b)\}},
\]
where $a,b\in(0,1)$ with $a<b$. Then $Z^\phi_\eta$ denotes the number of blocks in the fragmentation process whose size is a value in $(\nicefrac{\eta}{b},\nicefrac{\eta}{a})$. Let us remark that for general branching processes the respective characteristic gives the number of particles born in the time interval $(-\ln(\eta)+\ln(a),-\ln(\eta)+\ln(b))$, which is a finite number. However, if the fragmentation process has an infinite dislocation measure, we have $Z^\phi_\eta\in\{0,\infty\}$ $\mathbb P$-a.s. for any $a,b\in(0,1)$, so this particular characteristic is not of interest for us. In Section~\ref{ss.erc} we will present  more complicated examples which fit the purpose of this paper. The crucial point  when considering random characteristics for fragmentation processes with an infinite dislocation measure is that, in contrast to the situation of a finite dislocation measure,  $\phi$ needs to depend on $\pi$ as we need to integrate it with respect to the dislocation measure. 

Before we immerse ourselves in the proof of Theorem~\ref{t.2.1}, we first consider two applications of  this result.

\section{Examples}\label{ss.erc}
As a motivation for the main result of this paper let us now discuss two special random characteristics for which Theorem~\ref{t.2.1} is applicable. In the spirit of \cite{84} this section is concerned  with the energy that is needed to crush block in the mining industry, a theme that was taken up also in \cite{FKM10}. Moreover, we derive some asymptotic properties of the empirical mean associated with the stopped fragmentation process that we defined in Section~\ref{s.hfp}.

Recall that we consider  a self-similar fragmentation process $\Pi$ satisfying (\ref{e.cE}) as well as Hypothesis~\ref{h.1}. 

The main goal of this section is to prove almost sure convergence for the energy functional considered in \cite{84}. In that paper Bertoin and Mart\'inez proved limit theorems for the energy, properly discounted, that is required to fragment blocks of unit size to fragments of size less than some given value $\eta\in(0,1]$. In the setting of a finite dislocation measure they apply Corollary~3.3 of \cite{80} to prove $\mathscr L^1$-convergenc of this discounted energy functional as the value $\eta$ tends to zero. They further remark that Theorem~5.4 of \cite{80} can be used to obtain almost sure convergence under an appropriate additional assumption on the dislocation. However, the more interesting case is that of an infinite dislocation measure and in this spirit the main part of \cite{84} deals with an extension of the $\mathscr L^1$-convergence result to fragmentation processes with an infinite dislocation measure. Using some discretisation arguments Bertoin and Mart\'inez manage to prove $\mathscr L^1$-convergence, cf. Theorem~2 in \cite{84}, under the additional assumption that a parameter that appears in the theorem is greater than $\underline p$. In addition they prove $\mathscr L^2$-convergence, cf. Theorem~3 in \cite{84}, without that requirement on the parameter and without using a discretisation method but assuming that the fragmentation process is conservative. Here we prove almost sure convergence for the same functional as in \cite{84}, neither being restricted to parameters above $\underline p$ nor being restricted to the conservative setting. 
For this purpose, consider a measurable random function $\psi:\mathcal P\times\Omega\to\R$ that is independent of $\Pi$. In \cite{84} Bertoin and Mart\'inez interpret such a function  as the cost function that measures the cost of the energy that is used in an instantaneous dislocation. Having this function at hand, let us define for any $p<p^*$ a random function $\mathcal E_p:(0,1]\times\Omega\to\R$ by
\[
\mathcal E_p(\eta)=\sum_{t\in\mathcal I}\mathds1_{\{|\Pi_{\kappa(t_i)}(t_i-)|\ge \eta\}}|\Pi_{\kappa(t_i)}(t_i-)|^{1+p}\psi(\pi(t_i))
\]
for every $\eta\in(0,1]$. In the context of \cite{84} this function represents the total energy that is used in the process of crushing some block of unit size to blocks of sizes less than $\eta$. For a more precise description of this interpretation we refer to the explanations provided in \cite{84}. Note that the process stops when all blocks are of size less than $\eta$, since a block that is less than $\eta$  ``falls through a sieve'' and does not fragment further. 

For comparison let us state the result by Bertoin and Mart\'inez that we aim at extending in particular to almost sure convergence.
\begin{proposition}[Theorem~2 of \cite{84}]\label{t.BM05}
Assume that $\Pi$ is homogenous and let $p\in(\underline p,p^*)$. Further, let $\varphi:\mathcal S\to\R$ be a measurable function, with $\varphi((1,0,\ldots))=0$, that satisfies
\begin{equation}\label{e.nu-norm}
\int_{\mathcal S}|\varphi({\bf s})|\nu(\dd{\bf s})<\infty.
\end{equation} 
Then
\[
\eta^{p^*-p}\mathcal E_p(\eta)\to\frac{M_\infty(p^*)}{(p^*-p)\Phi'(p^*)}\int_{\mathcal S}\varphi({\bf s})\nu(\dd{\bf s})
\]
in $\mathscr L^1(\mathbb P)$ as $\eta\downarrow0$.
\end{proposition}

The following strong law of large numbers extends Proposition~\ref{t.BM05}.
\begin{theorem}\label{p.fragenergy}
Assume that 
\begin{equation}\label{e.mu-norm}
\int_\mathcal P\mathbb E(\psi(\pi))\mu(\dd\pi)<\infty
\end{equation}
and let $p<p^*$. 
Then
\begin{equation}\label{e.frag_energy}
\eta^{p^*-p}\mathcal E_p(\eta)\to\frac{\Lambda_0(p^*)}{\Phi'(p^*)(p^*-p)}\int_\mathcal P\mathbb E(\psi(\pi))\mu(\dd\pi)
\end{equation}
$\mathbb P$-a.s. and in $\mathscr L^1(\mathbb P)$ as $\eta\downarrow0$. 
\end{theorem}


By considering the special case that $\psi$ is only concerned with the asymptotic frequencies of $\pi\in\mathcal P$, i.e. $\psi$ depends on $\pi$ only via the sizes of its blocks, it follows from equation~(3) of \cite{HKK10}  that the assumption~(\ref{e.mu-norm}) in Theorem~\ref{p.fragenergy} does in particular cover (\ref{e.nu-norm}) of Proposition~\ref{t.BM05}.

As we will see in the proof of Theorem~\ref{p.fragenergy}, the left-hand side of (\ref{e.frag_energy}) can be considered as the process counted with a particular random characteristic.

Theorem~\ref{p.fragenergy} leads to the following result regarding the asymptotic behaviour of an empirical measure associated with the stopping line that determines the stopped process $(\lambda_\eta)_{\eta\in(0,1]}$. In this regard we also refer to Theorem~1 in \cite{HKK10} as well as Corollary~2  in \cite{84}. For related  considerations in the context of fragmentation chains, see Corollary~1 in \cite{84} and Theorem~3.1 in \cite{HK11}. 
\begin{theorem}\label{c.p.fragenergy}
Let $f$ be a random  measurable $\R$-valued function on $[0,1]$ which is independent of $\Pi$ and satisfies 
\begin{equation}\label{e.bounded_L1}
\sup_{x\in[0,1]}f(x)\in\mathscr L^1(\mathbb P).
\end{equation} 
Then 
\begin{equation}\label{e.convergence_empirical_mean}
\sum_{k\in\N}\lambda^{1+p^*}_{\eta,k}f\left(\frac{\lambda_{\eta,k}}{\eta}\right)\to \frac{\Lambda_0(p^*)}{\Phi'(p^*)}\int_{(0,1)}\mathbb E(f(u))\left(\int_\mathcal S\sum_{k\in\N}\mathds1_{\{s_k<u\}}s_k^{1+p^*}\nu(\dd{\bf s})\right)\frac{\dd u}{u}
\end{equation}
$\mathbb P$-a.s. in $\mathscr L^1(\mathbb P)$ as $\eta\downarrow0$.  
\end{theorem}


Note that 
\[
\sum_{k\in\N}\lambda^{1+p^*}_{\eta,k}f\left(\frac{\lambda_{\eta,k}}{\eta}\right)=\int_{(0,1)}f\dd\rho_\eta
\]
is the empirical mean of $f$ with respect to the empirical measure 
\[
\rho_\eta:=\sum_{k\in\N}\lambda^{1+p^*}_{\eta,k}\delta_{\frac{\lambda_{\eta,k}}{\eta}},
\]
where $\delta_x$ denotes the Dirac measure on $x\in(0,1)$. Further, observe that 
\[
\sum_{k\in\N}\lambda^{1+p^*}_{\eta,k}f\left(\frac{\lambda_{\eta,k}}{\eta}\right) = \eta^{1+p^*}\sum_{i\in\mathcal I}\phi^{(i)}\left(\frac{\eta}{|\Pi_{\kappa(t_i)}(t_i-)|},\pi(t_i)\right)=\eta^{1+p^*}Z_\eta^{\phi},
\]
where $f$ is a bounded and measurable $\rpn$-valued function on $[0,1]$ and the $\phi^{(i)}$ are independent copies of the random characteristic $\phi$ given by
\[
\phi\left(x,\pi\right)= \sum_{k\in\N}\mathds 1_{\left\{|\pi_k|<x\le1\right\}}x^{-(1+p^*)}|\pi_k|^{1+p^*}f\left(\frac{|\pi_k|}{x}\right)
\]
for every $x\in\rpn$ and $\pi\in\mathcal P$. Consequently, Theorem~\ref{c.p.fragenergy} fits into the general framework of the present paper. Moreover, as pointed out in \cite{84}, Theorem~\ref{c.p.fragenergy} is also related to the fragmentation energy model as it corresponds to the potential energy, where the cost function $\psi$ is given by $\psi(\pi)=\sum_{j\in\N}|\pi_j|^{1+ p}-1$ for some $p<1$ and all $\pi\in\mathcal P$. This will be explained in the proof of Theorem~\ref{c.p.fragenergy}.

{\bf Proof of Theorem~\ref{p.fragenergy}}
The idea is to consider an appropriate random characteristic and to apply Theorem~\ref{t.2.1}. For the time being, assume that $\alpha=0$, i.e. $\Pi$ is homogenous. Consider the  random characteristic $\phi$ given by
\begin{equation}\label{e.example.1}
\phi(x,\pi)=\mathds1_{\{x\in(0,1]\}}x^{-(1+p)}\psi\left(\pi\right)
\end{equation}
for all $x\in\rpn$ and $\pi\in\mathcal P$.
Then we have
\begin{align}\label{e.example.2}
\eta^{1+p^*}Z^\phi_\eta &= \eta^{1+p^*}\sum_{i\in\mathcal I}\phi^{(i)}\left(\frac{\eta}{|\Pi_{\kappa(t_i)}(t_i-)|},\,\pi(t_i)\right)\notag
\\[0.5ex]
&= \eta^{p^*-p}\sum_{i\in\mathcal I}\mathds1_{\{|\Pi_{\kappa(t_i)}(t_i-)|\ge \eta\}}|\Pi_{\kappa(t_i)}(t_i-)|^{1+p}\psi^{(i)}(\pi(t_i))
\\[0.5ex]
&= \eta^{p^*-p}\mathcal E_p(\eta)\notag
\end{align}
for every $\eta\in(0,1]$, where the $\phi^{(i)}$ are independent copies of $\phi$. 
Hence, since 
\[
\int_\mathcal P\mathbb E\left(\sup_{\eta\in(0,1]}\eta^{1+\tilde p}\phi\left(\eta,\pi\right)\right)\mu(\dd\pi)=\int_\mathcal P\mathbb E\left(\sup_{\eta\in(0,1]}\eta^{\tilde p-p}\psi\left(\pi\right)\right)\mu(\dd\pi)\le\int_{\mathcal P}\mathbb E(\psi(\pi))\mu(\dd\pi)<\infty
\]
for every $\tilde p\in(p,p^*)$, the characteristic $\phi$ given by (\ref{e.example.1}) satisfies (\ref{e.assumption_1_strong}) and thus we deduce from (\ref{e.example.2}) and  Theorem~\ref{t.2.1} that 
\begin{align*}
\lim_{\eta\downarrow0}\left(\eta^{p^*-p}\mathcal E_p(\eta)\right) &= \lim_{\eta\downarrow0}\left(\eta^{1+p^*}Z^\phi_\eta\right) 
\\[0.5ex]
&= \frac{\Lambda_0(p^*)}{\Phi'(p^*)}\int_{(0,1)}\mathbb E\left(\sum_{k\in\N}|\Pi_{k}(t)|^{1+p^*}\int_\mathcal P\int_{(0,1]}\rho^{p^*}\phi\left(\rho,\pi\right)\dd\rho\mu(\dd\pi)\right)\dd t
\\[0.5ex]
&= \frac{\Lambda_0(p^*)}{\Phi'(p^*)}\int_{(0,1)}\mathbb E\left(M_t(p^*)\right)\dd t\int_{(0,1]}\rho^{p^*-p-1}\dd\rho\int_\mathcal P\psi(\pi)\mu(\dd\pi) 
\\[0.5ex]
&= \frac{\Lambda_0(p^*)}{\Phi'(p^*)(p^*-p)}\int_\mathcal P\psi(\pi)\mu(\dd\pi)
\end{align*}
$\mathbb P$-a.s. and in $\mathscr L^1(\mathbb P)$, where in the final equality we have used that $M(p^*)$ is a unit-mean martingale.

The extension to the self-similar case with $\alpha\ne0$   follows analogously to Part~IV of the proof of Theorem~\ref{t.2.1}.
\hfill$\square$

The proof of Theorem~\ref{c.p.fragenergy} is based on Theorem~\ref{p.fragenergy} with a particular choice of $\psi$. 

{\bf Proof of Theorem~\ref{c.p.fragenergy}}
Our approach follows the lines of the argument outlined on page~569 of \cite{84}. We aim at showing the asserted convergence  for the case that the function $f$ is a certain kind of random power function and then some approximation arguments can be used to obtain the results for the case that $f$ satisfies the conditions of Theorem~\ref{c.p.fragenergy}. 

Consider first a random function $\mathfrak f:(0,1)\times\Omega\to(1,\infty)$ such that 
\[
\mathfrak f(x)=x^{\mathfrak p-p^*}
\] 
holds $\mathbb P$-a.s. for some random $\mathfrak p:\Omega\to(\underline p,p^*)$, which is independent of $\Pi$, and all $x\in(0,1)$. Note that 
\begin{align*}
& \sum_{l\in\N}\lambda^{1+p^*}_{\eta,l}\mathfrak f\left(\frac{\lambda_{\eta,l}}{\eta}\right) 
\\[0.5ex]
&= \eta^{p^*-\mathfrak p}\sum_{l\in\N}\lambda^{1+\mathfrak p}_{\eta,l}
\\[0.5ex]
&= \eta^{p^*-\mathfrak p}\sum_{t\in\mathcal I}\mathds1_{\{|\Pi_{\kappa(t_i)}(t_i-)|\ge \eta\}}|\Pi_{\kappa(t_i)}(t_i-)|^{1+\mathfrak p}\sum_{j\in\N}|\pi_j(t_i)|^{1+\mathfrak p}\mathds1_{\{|\Pi_{\kappa(t_i)}(t_i-)\cap\pi_j(t_i)|<\eta\}} 
\\[0.5ex]
&= \eta^{p^*-\mathfrak p}\sum_{t\in\mathcal I}\mathds1_{\{|\Pi_{\kappa(t_i)}(t_i-)|\ge \eta\}}|\Pi_{\kappa(t_i)}(t_i-)|^{1+\mathfrak p}\left(\sum_{j\in\N}|\pi_j(t_i)|^{1+\mathfrak p}-1\right)
\\[0.5ex]
&= \eta^{p^*-\mathfrak p}\mathcal E_{\mathfrak p}(\eta),
\end{align*}
where the cost function $\psi$, that appears in the definition of $\mathcal E_{\mathfrak p}$, is given by $\psi(\pi)=\sum_{j\in\N}|\pi_j|^{1+\mathfrak p}-1$ for every $\pi\in\mathcal P$. Observe that $\mathfrak p>\underline p$ implies that $\psi$ satisfies (\ref{e.mu-norm}) and as a consequence of Theorem~\ref{p.fragenergy} we thus obtain that
\[
\lim_{\eta\downarrow0}\sum_{l\in\N}\lambda^{1+p^*}_{\eta,l}\mathfrak f\left(\frac{\lambda_{\eta,l}}{\eta}\right)
\]
exists $\mathbb P$-a.s. and in $\mathscr L^1(\mathbb P)$. Therefore, since $\mathfrak p-p^*>-1$, and hence $\mathfrak f$ is $\mathbb P$-a.s. Lebesgue integrable, it follows from Lemma~3.4 of \cite{Kno11}  that 
\begin{equation}\label{e.f-power_fct}
\sum_{k\in\N}\lambda^{1+p^*}_{\eta,k}\mathfrak f\left(\frac{\lambda_{\eta,k}}{\eta}\right)\to \frac{\Lambda_0(p^*)}{\Phi'(p^*)}\int_{(0,1)}\mathbb E(\mathfrak f(u))\left(\int_\mathcal S\sum_{k\in\N}\mathds1_{\{s_k<u\}}s_k^{1+p^*}\nu(\dd{\bf s})\right)\frac{\dd u}{u}
\end{equation}
$\mathbb P$-a.s. and in $\mathscr L^1(\mathbb P)$ 
as $\eta\downarrow0$. Below we shall first prove the almost sure convergence and eventually deduce from this by means of the DCT  the $\mathscr L^1$-convergence. However, notice that in order to infer the above almost sure convergence we used having $\mathscr L^1$-convergence in Theorem~\ref{p.fragenergy}. 

Note that the martingale property of $\Lambda(p^*)$ implies that $\sum_{l\in\N}\lambda^{1+p^*}_{\eta,l} f\left(\frac{\lambda_{\eta,l}}{\eta}\right)\to c\Lambda_0(p^*)$ $\mathbb P$-a.s. as $\eta\downarrow$ if $f\equiv c\in\R$ is constant.  Consequently, we deduce from  (\ref{e.f-power_fct})  and an approximation of $f$ in the light of the Stone-Weierstra{\ss} theorem that the convergence in (\ref{e.convergence_empirical_mean}) holds $\mathbb P$-a.s. for each random, independent of $\Pi$, continuous function  $f$ with compact support  $[a,b]\subsetneq(0,1)$ that satisfies (\ref{e.bounded_L1}). Indeed, for any such function $f$ let $(g_k)_{k\in\N}$ be a sequence of  functions $g_k:[0,1]\times\Omega\to\R$ for which the convergence in (\ref{e.convergence_empirical_mean}) holds $\mathbb P$-a.s. and such that $g_k\to f$  uniformly on $[a,b]$ $\mathbb P$-a.s. as $k\to\infty$. For any $\epsilon>0$ there does then exist a $k_\epsilon:\Omega\to\N$ such that
\[
\sum_{k\in\N}\lambda^{1+p^*}_{\eta,k}g_k\left(\frac{\lambda_{\eta,k}}{\eta}\right)-\epsilon\Lambda_\eta(p^*)\le\sum_{k\in\N}\lambda^{1+p^*}_{\eta,k}f\left(\frac{\lambda_{\eta,k}}{\eta}\right)\le\sum_{k\in\N}\lambda^{1+p^*}_{\eta,k}g_k\left(\frac{\lambda_{\eta,k}}{\eta}\right)+\epsilon\Lambda_\eta(p^*)
\]
$\mathbb P$-a.s. for all $k\ge k_\epsilon$. Letting first $\eta\downarrow0$ and then $\epsilon\downarrow0$ results in
\begin{align*}
\lim_{\eta\downarrow0}\sum_{k\in\N}\lambda^{1+p^*}_{\eta,k} f\left(\frac{\lambda_{\eta,k}}{\eta}\right) &= \frac{\Lambda_0(p^*)}{\Phi'(p^*)}\lim_{k\to\infty}\int_{(0,1)}\mathbb E( g_k(u))\left(\int_\mathcal S\sum_{k\in\N}\mathds1_{\{s_k<u\}}s_k^{1+p^*}\nu(\dd{\bf s})\right)\frac{\dd u}{u}
\\[0.5ex]
&=  \frac{\Lambda_0(p^*)}{\Phi'(p^*)}\int_{(0,1)}\mathbb E( f(u))\left(\int_\mathcal S\sum_{k\in\N}\mathds1_{\{s_k<u\}}s_k^{1+p^*}\nu(\dd{\bf s})\right)\frac{\dd u}{u}
\end{align*}
$\mathbb P$-a.s., where in view of (\ref{e.bounded_L1}) the final equality follows from the DCT.
Resorting to the arguments  in Part~II of the proof of Theorem~3.2 in \cite{Kno11} we then conclude that the convergence in (\ref{e.convergence_empirical_mean}) holds $\mathbb P$-a.s. even  for every random  bounded and measurable function $f:[0,1]\times\Omega\to\R$ that is independent of $\Pi$ and satisfies (\ref{e.bounded_L1}). 

Since by means of Proposition~3.5 in \cite{Kno11} (cf. also Theorem~2 in \cite{89} for the conservative case) we have $\sup_{\eta\in[0,1]}\Lambda_\eta(p^*)\in\mathscr L^1(\mathbb P)$, it follows that
\[
\sup_{\eta\in[0,1]}\left(\sum_{k\in\N}\lambda^{1+p^*}_{\eta,k}f\left(\frac{\lambda_{\eta,k}}{\eta}\right)\right)\le\sup_{x\in[0,1]}f(x)\sup_{\eta\in[0,1]}\Lambda_\eta(p^*)\in\mathscr L^1(\mathbb P)
\]
and in view of the above shown almost sure convergence we thus deduce from the DCT that the convergence in (\ref{e.convergence_empirical_mean}) also holds in $\mathscr L^1(\mathbb P)$.
\hfill$\square$

\section{Proof of Theorem~\ref{t.2.1}}\label{ss.p.{t.2.1}}

This section is devoted to the proof of our main result. We first establish some auxiliary results to which we shall then resort in the proof of Theorem~\ref{t.2.1}.

Throughout this section, unless otherwise stated, assume that  $\Pi$ is homogenous. We will first prove the result in the homogenous case and then extend it to the general self-similar setting. Moreover, we assume that $\phi$ is nonnegative as the generalisation to an $\R$-valued $\phi$ is easily obtained by considering the positive and negative parts of $\phi$ separately and adding the two parts together.

In order to state the first auxiliary result we need to introduce some notation. To this end, for every $a>0$ let $\mathcal I^{(a)}$ be given by
\[
\mathcal I^{(a)}:=\left\{i\in\mathcal I:t_i\in(na)_{n\in\N}\right\}.
\]
In addition, set 
\[
\mathcal J^{(a)}_\eta:=
\{i\in\mathcal I^{(a)}:|\Pi_{\kappa(t_i)}(t_i-)|\ge \eta\}
\] 
as well as 
\[
T^{(a)}_\eta:=\text{card}\left(\mathcal J^{(a)}_\eta\right)
\] 
for each $a>0$ and $\eta\in(0,1]$. Further, adopt 
\[
T_{\eta,\rho}:=\text{card}\left(\left\{k\in\N:\lambda_{\eta,k}\ge \eta\rho\right\}\right)
\] 
for any $\eta,\rho\in(0,1]$. 

The following result is crucial for our further considerations.

\begin{proposition}\label{l.t.SLLN.fprc.1}
Let $a>0$. Then we have 
\begin{equation}\label{e.l.t.SLLN.fprc.1.statement}
\sup_{\eta\in(0,1]}\left(\eta^{1+p^*}T^{(a)}_\eta\right)\in\mathscr L^1(\mathbb P).
\end{equation}
\end{proposition}

\begin{pf}
The proof is divided into three parts. In the first part we show that 
\[
\limsup_{\eta\downarrow0}\left(\eta^{1+p^*}T^{(a)}_\eta\right) \in \mathscr L^1(\mathbb  P)
\]
and in the second part we prove that
\[
\zeta^{1+p^*}T^{(a)}_\zeta\in\mathscr L^1(\mathbb P)
\]
for every finite and nonnegative random variable $\zeta$. In the third part we combine these two conclusions in order to deduce that (\ref{e.l.t.SLLN.fprc.1.statement}) holds.

\underline{Part I}  
Throughout this proof we shall consider the random $\rpn$-valued function $\xi_a$ on $(0,1)$, given by 
\[
\xi_a(\rho):=\sum_{k\in\N}\mathds1_{\{|\pi_k(\tau_\rho)|\ge\rho\}}\mathds1_{\{\tau_\rho<a\}}+\mathds1_{\{\tau_\rho\ge a\}}
\]
for every $\rho\in(0,1)$. First, observe that by means of (\ref{e.cE}) and Proposition~2 in Section~0.5 of \cite{Ber96} we have 
\begin{align*}
& \lim_{\rho\downarrow0}\mathbb E\left(\frac{1}{\xi_a(\rho)+\epsilon}\right) 
\\[0.5ex]
&\le
 \lim_{\rho\downarrow0}\left(\frac{1}{2+\epsilon}\mathbb P\left(\xi_a(\rho)\ge2\right)+\frac{1}{1+\epsilon}\mathbb P\left(\xi_a(\rho)=1\right)+\frac{1}{\epsilon}\mathbb P\left(\xi_a(\rho)=0\right)\right)
\\[0.5ex]
&=
\lim_{\rho\downarrow0}\left(\frac{1}{2+\epsilon}\frac{\nu(s_2\ge\rho)}{\nu(s_1\le1-\rho)}
+\frac{1}{1+\epsilon}\frac{\nu(s_1\in[\rho,1-\rho],s_2<\rho)}{\nu(s_1\le1-\rho)}
+\frac{1}{\epsilon}\frac{\nu(s_1<\rho)}{\nu(s_1\le1-\rho)}\right)
\\[0.5ex]
&\qquad \cdot\lim_{\rho\downarrow0}\mathbb P({\bf e}_\rho<a)+\frac{1}{1+\epsilon}\lim_{\rho\downarrow0}\mathbb P({\bf e}_\rho\ge a)
\\[0.5ex]
&= \frac{1}{2+\epsilon}
\\[0.5ex]
&<\frac{1}{1+\epsilon},
\end{align*}
where  the random variable ${\bf e}_\rho$ is exponentially distributed with parameter $\nu({\bf s}\in\mathcal S:s_1\le1-\rho)$. Moreover,
\[
\lim_{\rho\downarrow0}\mathbb E(\xi_a(\rho)) \ge \lim_{\rho\downarrow0}\left(2\frac{\nu(s_2\ge\rho)}{\nu(s_1\le1-\rho)}
+\frac{\nu(s_1\in[\rho,1-\rho],s_2<\rho)}{\nu(s_1\le1-\rho)}\right)=2>1.
\]
In the light of the above estimates, for the remainder of this proof we fix some $\rho\in(0,1)$ such that
\begin{equation}\label{e.l.t.SLLN.fprc.1.0.b}
\mathbb E(\xi_a(\rho))>1\qquad\text{and}\qquad\mathbb E\left(\frac{1}{\xi_a(\rho)+\epsilon}\right)<\frac{1}{1+\epsilon}
\end{equation}
hold for every $\epsilon>0$. 

Observe that
\[
T_{\eta,\rho}\ge\sum_{j\in\mathcal J^{(a)}_\eta}\xi^{(j)}_a(\rho)-T^{(a)}_\eta
\]
$\mathbb P$-a.s., where the $\xi^{(j)}_a$ are independent copies of $\xi_a$. By means of Kolmogorov's strong law of large numbers this estimate results in 
\begin{equation}\label{e.l.t.SLLN.fprc.1.0}
\liminf_{\eta\downarrow0}\frac{T_{\eta,\rho}}{T^{(a)}_\eta}\ge\mathbb E\left(\xi_a(\rho)\right)-1>0
\end{equation}
$\mathbb P$-a.s. on $\{T^{(a)}_\eta\to\infty\text{ as }\eta\downarrow0\}$, where the positivity holds in view of (\ref{e.l.t.SLLN.fprc.1.0.b}). In the light of the estimate $T_{\eta,\rho}\le (\eta\rho)^{-(1+p^*)}\Lambda_\eta(p^*)$ $\mathbb P$-a.s., we thus infer that 
\begin{eqnarray}\label{e.l.t.SLLN.fprc.1.1}
\limsup_{\eta\downarrow0}\left(\eta^{1+p^*}T^{(a)}_\eta\right) & \le & \liminf_{\eta\downarrow0}\left(\eta^{1+p^*}T_{\eta,\rho}\right)(\mathbb E(\xi_a(\rho))-1)^{-1}\notag
\\[0.5ex]
& \le & {\rho}^{-(1+p^*)}\left(\mathbb E(\xi_a(\rho))-1\right)^{-1}\Lambda_0(p^*)
\end{eqnarray}
$\mathbb P$-almost surely. 

\underline{Part II} For the remainder of this proof fix some $\epsilon>0$. Further, throughout this second part of the proof 
let $\zeta:\Omega\to(0,1]$ be some finite and nonnegative random variable. Along the lines leading to (\ref{e.l.t.SLLN.fprc.1.1}) we obtain that

\begin{align}\label{e.l.t.SLLN.fprc.1.1b}
{\rho}^{-(1+p^*)}\Lambda_\zeta(p^*) &\ge \zeta^{1+p^*}\sum_{j\in\mathcal J^{(a)}_\zeta}\xi^{(j)}_a(\rho)-\zeta^{1+p^*}T^{(a)}_\zeta\notag
\\[0.5ex]
&= \zeta^{1+p^*}\sum_{j\in\mathcal J^{(a)}_\zeta}\left(\xi^{(j)}_a(\rho)+\epsilon\right)-(1+\epsilon)\zeta^{1+p^*}T^{(a)}_\zeta
\end{align}
$\mathbb P$-almost surely. Let $q\in[\nicefrac{1}{2},1)$ and let $p:=\nicefrac{q}{(q-1)}\in[-1,0)$ be its conjugate, that is to say $p^{-1}+q^{-1}=1$. Then we infer that
\begin{align}\label{e.l.t.SLLN.fprc.1.1b2}
& \mathbb E\left(\zeta^{1+p^*}\sum_{i\in\mathcal J^{(a)}_\zeta}\left(\xi^{(j)}_a(\rho)+\epsilon\right)\right)\notag
\\[0.5ex]
&= \sum_{i\in\N}\mathbb E\left(\mathds1_{\{i\in\mathcal J^{(a)}_\zeta\}}\zeta^{1+p^*}\left(\xi^{(j)}_a(\rho)+\epsilon\right)\right)\notag
\\[0.5ex]
&\ge \sum_{i\in\N}\mathbb E\left(\left[\mathds1_{\{i\in\mathcal J^{(a)}_\zeta\}}\zeta^{1+p^*}\right]^q\right)^{\nicefrac{1}{q}}\mathbb E\left(\left[\xi^{(j)}_a(\rho)+\epsilon\right]^p\right)^{\nicefrac{1}{p}}\notag
\\[0.5ex]
&\ge \mathbb E\left(\left[\xi_a(\rho)+\epsilon\right]^p\right)^{\nicefrac{1}{p}}
\sum_{i\in\N}\mathbb E\left(\mathds1_{\{i\in\mathcal J^{(a)}_\zeta\}}\zeta^{1+p^*}\right)\mathbb E\left(\mathds1_{\{i\in\mathcal J^{(a)}_\zeta\}}\zeta^{1+p^*}\right)^{\frac{1}{q}-1}
\\[0.5ex]
&\ge \mathbb E\left(\left[\xi_a(\rho)+\epsilon\right]^p\right)^{\nicefrac{1}{p}}
\left[\sum_{i\in\N}\mathbb E\left(\mathds1_{\{i\in\mathcal J^{(a)}_\zeta\}}\zeta^{1+p^*}\right)\right]^{2-\frac{1}{q}}\notag
\left[\sum_{i\in\N}\mathbb E\left(\mathds1_{\{i\in\mathcal J^{(a)}_\zeta\}}\zeta^{1+p^*}\right)\right]^{\frac{1}{q}-1}\notag
\\[0.5ex]
&= \mathbb E\left(\left[\xi_a(\rho)+\epsilon\right]^p\right)^{\nicefrac{1}{p}}\mathbb E\left(\zeta^{1+p^*}\sum_{i\in\N}\mathds1_{\{i\in\mathcal J^{(a)}_\zeta\}}\right)\notag
\\[0.5ex]
&= \mathbb E\left(\left[\xi_a(\rho)+\epsilon\right]^p\right)^{\nicefrac{1}{p}}\mathbb E\left(\zeta^{1+p^*}T^{(a)}_\zeta\right).\notag
\end{align}
Notice that in the first equality of (\ref{e.l.t.SLLN.fprc.1.1b2}) we used the MCT, in the subsequent step we resorted to the reverse H\"older inequality (see \cite{110}) and the last lower estimate is a consequence of Jensen's inequality. Let us mention that in view of $p<0$ we use the $\epsilon$ in order to avoid having ``$0^p$''. 

Observe that according to (\ref{e.l.t.SLLN.fprc.1.0.b}) and in the light  of $\xi_a(\rho)\le\rho^{-1}$ we have
\begin{equation}\label{e.finiteness}
\mathbb E\left(\left(\xi_a(\rho)+\epsilon\right)^{-1}\right)^{-1}-(1+\epsilon)\in(0,\infty).
\end{equation}

It follows from Proposition~3.5 in \cite{Kno11} (cf. also Theorem~2 in \cite{89} for the conservative case)  
that $\sup_{\eta\in(0,1]}\Lambda_\eta(p^*)\in\mathscr L^1(\mathbb P)$ and thus we infer from   (\ref{e.l.t.SLLN.fprc.1.1b}) and (\ref{e.l.t.SLLN.fprc.1.1b2}) that
\begin{align*}
{\rho}^{-(1+p^*)}\mathbb E\left(\sup_{\eta\in(0,1]}\Lambda_\eta(p^*)\right) &\ge \mathbb E\left({\rho}^{-(1+p^*)}\Lambda_\zeta(p^*)\right)
\\[0.5ex]
&\ge\mathbb E\left(\zeta^{1+p^*}T^{(a)}_\zeta\right)\left(\mathbb E\left(\left(\xi_a(\rho)+\epsilon\right)^{-1}\right)^{-1}-(1+\epsilon)\right),
\end{align*}
which by means of  (\ref{e.finiteness}) results in
\begin{equation}\label{e.l.t.SLLN.fprc.1.2}
\mathbb E\left(\zeta^{1+p^*}T^{(a)}_\zeta\right) \le  \rho^{-(1+p^*)}\mathbb E\left(\sup_{\eta\in(0,1]}\Lambda_\eta(p^*)\right)\left(\mathbb E\left(\left(\xi_a(\rho)+\epsilon\right)^{-1}\right)^{-1}-(1+\epsilon)\right)^{-1}<\infty.
\end{equation}

\underline{Part III}
Let $\delta>0$ and in view of (\ref{e.l.t.SLLN.fprc.1.1}) let $\zeta_\delta$ be a $(0,1]$-valued random variable such that 
\[
\zeta_\delta^{1+p^*}T^{(a)}_{\zeta_\delta}\ge\sup_{\eta\in(0,1]}\left(\eta^{1+p^*}T^{(a)}_\eta\right)-\delta
\]
on the event $\{\sup_{\eta\in(0,1]}(\eta^{1+p^*}T^{(a)}_\eta)\neq\limsup_{\eta\to0}(\eta^{1+p^*}T^{(a)}_\eta)\}$. Then, by means of (\ref{e.l.t.SLLN.fprc.1.1}) and (\ref{e.l.t.SLLN.fprc.1.2}), we obtain
\begin{align*}
& \mathbb E\left(\sup_{\eta\in(0,1]}\left(\eta^{1+p^*}T^{(a)}_\eta\right)\right) 
\\[0.5ex]
&\le \mathbb E\left(\limsup_{\eta\to0}\left(\eta^{1+p^*}T^{(a)}_\eta\right)\right)+\mathbb E\left(\zeta_\delta^{1+p^*}T^{(a)}_{\zeta_\delta}+\delta\right)
\\[0.5ex]
&\le {\rho}^{-(1+p^*)}(\mathbb E(\xi_a(\rho))-1)^{-1}+\rho^{-(1+p^*)}\mathbb E\left(\sup_{\eta\in(0,1]}\Lambda_\eta(p^*)\right)\left(\mathbb E\left(\left(\xi_a(\rho)+\epsilon\right)^{-1}\right)^{-1}-(1+\epsilon)\right)^{-1}+\delta
\\[0.5ex]
&< \infty,
\end{align*}
which shows that (\ref{e.l.t.SLLN.fprc.1.statement}) holds and thus completes the proof.
\end{pf}

Let us continue with establishing some auxiliary results.

\begin{lemma}\label{l.B}
Let $\phi$ be a random characteristic satisfying (\ref{e.assumption_1}). Then the limit $\lim_{\eta\downarrow0}\mathbb E(\eta^{1+p^*}Z^{\phi}_{\eta})$ exists and satisfies
\[
\lim_{\eta\downarrow0}\mathbb E\left(\eta^{1+p^*}Z^{\phi}_{\eta}\right)=\frac{1}{\Phi'(p^*)}\int_{(0,1)}\mathbb E\left(\sum_{k\in\N}|\Pi_{k}(s)|^{1+p^*}\int_\mathcal P\int_{(0,1]}\rho^{p^*}\phi\left(\rho,\pi\right)\dd\rho\mu(\dd\pi)\right)\dd s.
\]
\end{lemma}

\begin{pf}
Consider a function $g:\rpn\to\rpn$ defined by
\[
g(t):=\mathbb E\left(e^{-(1+p^*)t}Z^{\phi}_{e^{-t}}\right)
\]
for each $t\in\rpn$. Then we have
\begin{align*}
g(t)\notag
&= e^{-(1+p^*)t}\mathbb E\left(\sum_{i\in\mathcal I}\mathds1_{\{t_i\le 1\}}{\phi^{(i)}}\left(e^{-(t+\ln(|\Pi_{\kappa(t_i)}(t_i-)|))},\pi(t_i)\right)+\mathbb E\left(\left.\sum_{k\in\N}Z^{\phi,k}_{e^{-(t+\ln(|\Pi_k(1)|))}}\right|\mathscr F_1\right)\right)\notag
\\[0.5ex]
&= e^{-(1+p^*)t}\mathbb E\left(\sum_{i\in\mathcal I}\mathds1_{\{t_i\le 1\}}{\phi^{(i)}}\left(e^{-(t+\ln(|\Pi_{\kappa(t_i)}(t_i-)|))},\pi(t_i)\right)\right)
\\[0.5ex]
&\qquad +\mathbb E\left(\sum_{k\in\N}|\Pi_k(1)|^{1+p^*}g\left(t+\ln(|\Pi_k(1)|)\right)\right),\notag
\end{align*}
where the $\phi^{(i)}$ are independent copies of $\phi$ and, given $\mathscr F_1$, the $Z^{\phi,k}$ are independent copies of $Z^{\phi}$. 
That is, $g$  satisfies the renewal equation
\begin{equation}\label{e.renewal_equation.2}
g(t)=f(t)+(g*\varrho)(t) 
\end{equation}
for every $t\in\rpn$, where $f:\rpn\to\rpn$ is given by
\begin{equation}\label{e.renewal_equation.3}
f(t)=e^{-(1+p^*)t}\mathbb E\left(\sum_{i\in\mathcal I}\mathds1_{\{t_i\le 1\}}{\phi^{(i)}}\left(e^{-(t+\ln(|\Pi_{\kappa(t_i)}(t_i-)|))},\pi(t_i)\right)\right)
\end{equation}
for any $t\in\rpn$ and the measure $\varrho$ is defined by
\[
\varrho(\dd x) = \mathbb E\left(\sum_{k\in\N}|\Pi_k(1)|^{1+p^*}\mathds1_{\{-\ln(|\Pi_k(1)|)\in\dd x\}}\right)
\]
for every $x\in\rpn$. Note that
\[
\sum_{i\in\mathcal I}\mathds1_{\{t_i\le 1\}}{\phi^{(i)}}\left(e^{-(t+\ln(|\Pi_{\kappa(t_i)}(t_i-)|))},\pi(t_i)\right)\le\sum_{j\in\mathcal J^{(a)}_{\rho}}\sum_{i\in\mathcal I^{(j)}}\mathds1_{\{t_{i,j}\le a\}}\phi^{(i,j)}\left(\frac{\eta}{|\Pi_{\kappa(t_{i,j})}(t_{i,j}-)|},\pi(t_{i,j})\right)
\]
for all $a\ge1$ and $t\in\rp$, $\eta:=e^{-t}$ and any $\rho\in(\eta,1)$, where the $(t_{i,j})_{i\in\mathcal I^{(j)}}$ and $\phi^{(i,j)}$  are independent copies of $(t_i)_{i\in\mathcal I}$ and $\phi$ respectively. Therefore, we infer from the upcoming argument in (\ref{e.c.t.SLLN.1.3.0.a2}) that there exists some $\tilde p\in(\underline p,p^*)$ such that
\begin{equation}\label{e.renewal_equation.4}
f(t)=\mathcal O\left(e^{-(p^*-\tilde p)t}\right)
\end{equation}
as $t\to\infty$, where  $f$ is given by (\ref{e.renewal_equation.3}) and $\mathcal O$ is the Bachmann-Landau notation. Moreover,  the compensation formula for Poisson point processes, in conjunction with Tonelli's theorem, implies that
\begin{align}\label{e.renewal_equation.5}
\sup_{t\in[0,t_0)}f(t) &\le \sup_{t\in[0,t_0)}\mathbb E\left(\sum_{i\in\mathcal I}\mathds1_{\{t_i\le1\}}\sup_{s\in[0,t)}\left(e^{-(1+p^*)s}\phi^{(i)}\left(e^{-s},\pi(t_i)\right)\right)\right)\notag
\\[0.5ex]
&\le e^{\beta t_0}\int_\mathcal P\mathbb E\left(\sup_{s\in[0,t_0)}\left(e^{-(1+p^*+\beta)s}\phi\left(e^{-s},\pi\right)\right)\right)\mu(\dd\pi)
\\[0.5ex]
&<\infty\notag
\end{align}
for all $\beta>0$ and $t_0\in\rpn$, where the finiteness is a consequence of (\ref{e.assumption_1}).
Hence, since $t\mapsto e^{-(p^*-\tilde p)t}$  is directly Riemann integrable, we deduce from (\ref{e.renewal_equation.4}) and (\ref{e.renewal_equation.5}) that the function $f$ is directly Riemann integrable, see Lemma 3.4.1 in \cite{Asm03}. 

Observe that the compensation formula and Tonelli's theorem yield that
\begin{align*}
& \int_{\rpn}f(t)\dd t 
\\[0.5ex]
&= \mathbb E\left(\int_{(0,1)}\int_\mathcal P\int_{\rpn}e^{-(1+p^*)t}\phi\left(e^{-(t+\ln(|\Pi_{k(s)}(s)|))},\pi\right)\dd t\mu(\dd\pi)\dd s\right)
\\[0.5ex]
&= \int_{(0,1)}\mathbb E\left(\sum_{k\in\N}|\Pi_{k}(s)|^{1+p^*}\int_\mathcal P\int_{\rpn}e^{-(1+p^*)(t+\ln(|\Pi_{k}(s)|))}\phi\left(e^{-(t+\ln(|\Pi_{k}(s)|))},\pi\right)\dd t\mu(\dd\pi)\right)\dd s
\\[0.5ex]
&= \int_{(0,1)}\mathbb E\left(\sum_{k\in\N}|\Pi_{k}(s)|^{1+p^*}\int_\mathcal P\int_{\rpn}e^{-(1+p^*)u}\phi\left(e^{-u},\pi\right)\dd u\mu(\dd\pi)\right)\dd s
\end{align*}
In order to complete the proof recall that under $\mathbb P$ the process $(-\ln(|\Pi_1(t)|))_{t\in\rpn}$, is a subordinator with Laplace exponent $\Phi$ and consider the change of measure
\[
\left.\frac{\dd\mathbb P^{(p^*)}}{\dd\mathbb P}\right|_{\mathscr F_t}=|\Pi_1(t)|^{p^*}.
\]
Let $\mathbb E^{(p^*)}$ denote the expectation under $\mathbb P^{(p^*)}$. Below we resort to the following many-to-one identity for fragmentations (see Lemma~2 in \cite{HKK10}):
\begin{equation}\label{e.many-to-1}
\mathbb E\left(\sum_{k\in\N}|\Pi_k(t)|^{1+p^*}f(|\Pi_k(t)|)\right)=\mathbb E^{(p^*)}\left(f(|\Pi_k(t)|)\right)
\end{equation}
for every measurable $f:[0,1]\to\rpn$ and $t\in\rpn$. Since under the measure $\mathbb P^{(p^*)}$ the process $(-\ln(|\Pi_1(t)|))_{t\in\rpn}$ is a subordinator with Laplace exponent $\Phi_{p^*}(\lambda)=\Phi(\lambda+p^*)$, thus resulting in $\Phi_{p^*}'(0+)=\Phi'(p^*)$, we deduce from (\ref{e.many-to-1}) that
\begin{align*}
\int_{\rpn}x\varrho(\dd x) &=\int_{\rpn} \mathbb E\left(x\sum_{k\in\N}|\Pi_k(1)|^{1+p^*}\mathds1_{\{-\ln(|\Pi_k(1)|)\in\dd x\}}\right)
\\[0.5ex]
&= \mathbb E\left(|\Pi_1(1)|^{p^*}\int_{\rpn}x\mathds1_{\{-\ln(|\Pi_1(1)|)\in\dd x\}}\right)
\\[0.5ex]
&= \mathbb E^{(p^*)}\left(-\ln(|\Pi_1(1)|)\right)
\\[0.5ex]
&= \Phi'(p^*).
\end{align*}
In the light of (\ref{e.renewal_equation.2}) and the direct Riemann integrability of $f$ it follows from Smith's Key Renewal Theorem (cf. Theorem~A4.3 in \cite{EKM97}) that
\begin{align*}
\lim_{t\to\infty}\mathbb E\left(e^{-(1+p^*)t}Z^\phi_{e^{-t}}\right) &= \frac{\int_{\rpn}f(t)\dd t }{\int_{\rpn}x\varrho(\dd x)} 
\\[0.5ex]
&= \frac{1}{\Phi'(p^*)}\int_{(0,1)}\mathbb E\left(\sum_{k\in\N}|\Pi_{k}(s)|^{1+p^*}\int_\mathcal P\int_{\rpn}e^{-(1+p^*)u}\phi\left(e^{-u},\pi\right)\dd u\mu(\dd\pi)\right)\dd s
\end{align*}
and consequently
\[
\lim_{\eta\downarrow0}\mathbb E\left(\eta^{1+p^*}Z^{\phi}_{\eta}\right)=\frac{1}{\Phi'(p^*)}\int_{(0,1)}\mathbb E\left(\sum_{k\in\N}|\Pi_{k}(s)|^{1+p^*}\int_\mathcal P\int_{(0,1]}\rho^{1+p^*}\phi\left(\rho,\pi\right)\frac{\dd\rho}{\rho}\mu(\dd\pi)\right)\dd s.
\]
For further information regarding renewal theory and direct Riemann integrability we also refer to Section~3.10 in \cite{Res92}.
\end{pf}

For the remainder of this paper fix some $\iota>1$ and  for every random characteristic $\phi$  set
\begin{equation}\label{e.c.t.SLLN.1.2b.1}
\phi_{\iota,\eta,s}(x,\pi):=\phi(x,\pi)\mathds 1_{\{x>\eta^{s(\iota-1)}\}}\qquad\text{as well as}\qquad\phi_{\iota,\eta}:=\phi_{\iota,\eta,1}
\end{equation}
for all $s\ge1$, $x\in[0,1]$ and $\pi\in\mathcal P$. 

In the proof of Theorem~\ref{t.2.1} we shall make use of the following result: 
\begin{proposition}\label{t.appendix}
Let $\phi$ be a random characteristic satisfying (\ref{e.assumption_1}). Then we have
\[
\lim_{k\to\infty}\rho^{k\delta\iota(1+p^*)}Z^{\phi_{\iota,\rho^{k\delta}}}_{\rho^{k\delta\iota}}=\Lambda_0(p^*)\lim_{\eta\downarrow0}\mathbb E\left(\eta^{1+p^*}Z^\phi_\eta\right)
\]
$\mathbb P$-a.s. for all $\rho\in(0,1)$ and $\delta>0$.
\end{proposition}

We are now in a position to reap the fruits of our work and to prove the main result of this paper, resorting to Proposition~\ref{t.appendix} whose proof will be provided in Section~\ref{s.app}.

{\bf Proof of Theorem~\ref{t.2.1}}
Let $a>0$ and let $\tilde p\in(\underline p,p^*)$ be given by (\ref{e.assumption_1}). The proof is divided into four parts. In the first part of this proof we prove that the asserted almost sure convergence holds along ln-lattice times and in the second part the almost sure convergence along the real numbers in shown. Subsequently, in the third part we prove the $\mathscr L^1$-convergence. In these first three parts we wok under the assumption that $\Pi$ is homogenous. Finally, in the fourth part we show that we can drop this homogeneity assumption and consequently obtain the desired almost sure and $\mathscr L^1$-convergence for $\Pi$ being self-similar.

\underline{Part I} We start by proving the desired result along $\ln$-lattice times. 
To this end, set
\[
A_r:=\sup_{\eta\in(0,\rho^{(\iota-1)r}]}\int_\mathcal P\eta^{1+\tilde p}\mathbb E\left(\phi\left(\eta,\pi\right)\right)\mu(\dd\pi)
\]
for any $r\in\rp$ and observe that by means of the compensation formula for Poisson point processes we have
\begin{align*}
E^{a,j}_{r,k}(\rho) &:= \left.\mathbb E\left(\sum_{i\in\mathcal I}\mathds1_{\{t_i\le a\}}\mathds1_{\{|\Pi_{\kappa(t_i)}(t_i-)|\ge\rho^{r-k}\}}\phi^{(i)}\left(\frac{\rho^{\iota r}}{u_j|\Pi_{\kappa(t_i)}(t_i-)|},\pi(t_i)\right)\right)\right|_{u_j=|\Pi(\kappa(t_j))(t_j-)|}\notag
\\[0.5ex]
&=\int\limits_{(0,a)}\left.\mathbb E\left(\mathds1_{\{|\Pi_{\kappa(t)}(t-)|\ge\rho^{r-k}\}}\int_\mathcal P\phi\left(\frac{\rho^{\iota r}}{u_j|\Pi_{\kappa(t)}(t-)|},\pi\right)\mu(\dd\pi)\right)\right|_{u_j=|\Pi(\kappa(t_j))(t_j-)|}\dd s\notag
\\[0.5ex]
&= \left.\int\limits_{(0,a)}\int_\mathcal P\mathbb E\left(\mathds1_{\{|\Pi_{\kappa(t)}(t-)|\ge\rho^{r-k}\}}\left.\mathbb E\left(\phi\left(\frac{\rho^{\iota r}}{u_jv_t},\pi\right)\right)\right|_{v_t=|\Pi(\kappa(t))(t-)|}\right)\right|_{u_j=|\Pi(\kappa(t_j))(t_j-)|}\mu(\dd\pi)\dd s\notag
\\[0.5ex]
&\le aA_r\rho^{-(\iota r-k)(1+\tilde p)}
\end{align*}
for all $j,k\in\N$, $r\in\rp$  and $\rho\in(0,1)$. Hence, 
we deduce from the fragmentation property that
\begin{align*}
\mathbb E\left(\rho^{\iota r(1+p^*)}\left|Z^{\phi}_{\rho^{\iota r}}-Z^{\phi_{\iota,\rho^r}}_{\rho^{\iota r}}\right|\right)\notag
&=  \mathbb E\left(\sum_{j\in\mathcal J^{(a)}_{\rho^r}}\sum_{i\in\mathcal I^{(j)}}\mathds1_{\{t_{i,j}\le a\}}\rho^{\iota r(1+p^*)}\phi^{(i,j)}\left(\frac{\rho^{\iota r}}{|\Pi_{\kappa(t_{i,j})}(t_{i,j}-)|},\pi(t_{i,j})\right)\right)\notag
\\[0.5ex]
&= \sum_{k=0}^{r-1}\rho^{\iota r(1+p^*)}\mathbb E\left(\sum_{j\in\mathcal J^{(a)}_{\rho^{k+1}}\setminus\mathcal J^{(a)}_{\rho^k}}E^{a,j}_{r,k}(\rho)\right)\notag
\\[0.5ex]
&\le aA_r\sum_{k=0}^{r-1}\rho^{(\iota r-k)(p^*-\tilde p)}\mathbb E\left(\rho^{k(1+p^*)}\sharp\left(\mathcal J^{(a)}_{\rho^{k+1}}\setminus\mathcal J^{(a)}_{\rho^k}\right)\right)\notag
\\[0.5ex]
&\le aA_r\mathbb E\left(\sup_{u\in(0,1]}\left(u^{1+p^*}T^{(a)}_{u}\right)\right) \sum_{n=\left\lfloor(\iota-1)r+1\right\rfloor}^{\left\lceil\iota r\right\rceil}\rho^{n(p^*-\tilde p)}
\end{align*}
holds for every $r\in\rp$, where the $(t_{i,j})_{i\in\mathcal I^{(j)}}$ and $\phi^{(i,j)}$  are independent copies of $(t_i)_{i\in\mathcal I}$ and $\phi$ respectively.
Hence,  we infer that 
\begin{equation}\label{e.c.t.SLLN.1.3.0.a2}
\mathbb E\left(\rho^{\iota r(1+p^*)}\left|Z^{\phi}_{\rho^{\iota r}}-Z^{\phi_{\iota,\rho^r}}_{\rho^{\iota r}}\right|\right)
\le aA_r\rho^{-(1+p^*)}\mathbb E\left(\sup_{u\in(0,1]}\left(u^{1+p^*}T^{(a)}_{u}\right)\right)\frac{\rho^{(\iota-1)(p^*-\tilde p)r}}{1-\rho^{p^*-\tilde p}}
\end{equation}
for all $r\in\rp$.
For any $\delta,\epsilon>0$ let $n_{\delta,\epsilon}\in\N$ be such that 
\[
\sup_{\eta\in(0,\rho^{(\iota-1)\delta n}]}\int_\mathcal P\eta^{1+\tilde p}\mathbb E\left(\phi\left(\eta,\pi\right)\right)\mu(\dd\pi)\le\limsup_{\eta\downarrow0}\int_\mathcal P\eta^{1+\tilde p}\mathbb E\left(\phi\left(\eta,\pi\right)\right)\mu(\dd\pi)+\epsilon
\]
holds for all $n\ge n_{\delta,\epsilon}$. 
Resorting to the Chebyshev-Markov inequality we thus deduce in view of (\ref{e.c.t.SLLN.1.3.0.a2})  that 
\begin{align*}
& \sum_{n=n_{\delta,\epsilon}}^\infty\mathbb P\left(\rho^{\iota \delta n(1+p^*)}\left|Z^{\phi}_{\rho^{\iota \delta n}}-Z^{\phi_{\iota,\rho^{\delta n}}}_{\rho^{\iota {\delta n}}}\right|\ge\varepsilon\right) 
\\[0.5ex]
&\le \frac{1}{\varepsilon}\sum_{n=n_{\delta,\epsilon}}^\infty\mathbb E\left(\rho^{\iota \delta n(1+p^*)}\left|Z^{\phi}_{\rho^{\iota \delta n}}-Z^{\phi_{\iota,\rho^{\delta n}}}_{\rho^{\iota {\delta n}}}\right|\right)
\\[0.5ex]
&\le \frac{1}{\varepsilon}\rho^{-(1+p^*)}\mathbb E\left(\sup_{u\in(0,1]}\left(u^{1+p^*}T^{(a)}_{u}\right)\right)\frac{1}{1-\rho^{p^*-\tilde p}}\sum_{n=n_{\delta,\epsilon}}^\infty A^a_{\delta n}\rho^{(\iota-1)(p^*-\tilde p)\delta n}
\\[0.5ex]
&\le \frac{1}{\varepsilon}\rho^{-(1+p^*)}
\mathbb E\left(\sup_{u\in(0,1]}\left(u^{1+p^*}T^{(a)}_{u}\right)\right)\frac{1}{1-\rho^{p^*-\tilde p}}
\\[0.5ex]
&\qquad \cdot a\left(\limsup_{\eta\downarrow0}\int_\mathcal P\eta^{1+\tilde p}\mathbb E\left(\phi\left(\eta,\pi\right)\right)\mu(\dd\pi)+\epsilon\right)\sum_{n=n_{\delta,\epsilon}}^\infty\rho^{(\iota-1)(p^*-\tilde p)\delta n}
\\[0.5ex]
&< \infty
\end{align*}
for all $\delta,\epsilon,\varepsilon>0$, where the finiteness follows from (\ref{e.assumption_1}), Proposition~\ref{l.t.SLLN.fprc.1} and the fact that the geometric series $\sum_{n\in\N}\rho^{(\iota-1)(p^*-\tilde p)\delta n}$ is finite. Therefore,  the Borel-Cantelli lemma yields that
\[
\lim_{n\to\infty}\left(\rho^{\iota \delta n(1+p^*)}\left|Z^{\phi}_{\rho^{\iota \delta n}}-Z^{\phi_{\iota,\rho^{\delta n}}}_{\rho^{\iota {\delta n}}}\right|\right)=0
\]
$\mathbb P$-a.s. for every $\delta>0$. Consequently, by means of Lemma~\ref{l.B} and Proposition~\ref{t.appendix} we infer that
\begin{align}\label{e.c.t.SLLN.1.3.0}
& \lim_{k\to\infty}\rho^{\iota \delta k(1+p^*)}Z^\phi_{\rho^{\iota \delta k}}-\frac{\Lambda_0(p^*)}{\Phi'(p^*)}\int_{(0,1)}\mathbb E\left(\sum_{k\in\N}|\Pi_{k}(s)|^{1+p^*}\int_\mathcal P\int_{(0,1]}\rho^{p^*}\phi\left(\rho,\pi\right)\dd\rho\mu(\dd\pi)\right)\dd s\notag
\\[0.5ex]
&\le \lim_{k\to\infty}\left(\rho^{\iota \delta k(1+p^*)}\left|Z^\phi_{\rho^{\iota \delta k}}-Z^{\phi_{\iota,\rho^{\delta k}}}_{\rho^{\iota {\delta k}}}\right|\right)+\left|\lim_{k\to\infty}\rho^{\iota \delta k(1+p^*)}Z^{\phi_{\iota,\eta^{k\delta}}}_{\rho^{\iota \delta k}}-\Lambda_0(p^*)\lim_{\eta\downarrow0}\mathbb E\left(\eta^{1+p^*}Z^\phi_\eta\right)\right|\notag
\\[0.5ex]
& =0
\end{align}
holds $\mathbb P$-a.s. for all $\delta>0$ and $\rho\in(0,1)$.

\underline{Part II}
Let us now extend the convergence along $\ln$-lattice sequences in (\ref{e.c.t.SLLN.1.3.0}) to convergence along the positive real numbers. To this end, fix some $\rho\in(0,1)$ and  define $\underline\phi_\delta$ as well as $\bar\phi_\delta$ by
\[
\underline\phi_\delta\left(\eta,\pi\right)=\inf_{u\in(0,\delta)}\phi\left(\eta\rho^u,\pi\right)\quad\text{and}\quad\bar\phi_\delta\left(\eta,\pi\right)=\sup_{u\in(0,\delta)}\phi\left(\eta\rho^{-u},\pi\right)
\]
for any $\eta\in(0,1]$ and $\pi\in\mathcal P$. Notice that $\phi$ satisfying (\ref{e.assumption_1}) implies that also $\underline\phi_\delta$ and $\bar\phi_\delta$ satisfy (\ref{e.assumption_1}). For the time being, let $\delta>0$ and for any $r\in\rpn$   let $n_r\in\N$ be the unique natural number such that $r\in(n_r\delta,(n_r+1)\delta)$. Then we have
\begin{equation}\label{e.t.SLLN.fprc.II.1}
\rho^{(n_r+1)\delta(1+p^*)}Z^{\underline\phi_\delta}_{\rho^{(n_r+1)\delta}}\le \rho^{r(1+p^*)}Z^\phi_{\rho^r}\le \rho^{n_r\delta(1+p^*)}Z^{\bar\phi_\delta}_{\rho^{n_r\delta}}
\end{equation}
for all $r\in\rpn$. According to (\ref{e.c.t.SLLN.1.3.0}) we have, by Tonelli's theorem,
\[
\lim_{t\to\infty}\left(\rho^{(n_r+1)\delta(1+p^*)}Z^{\underline\phi_\delta}_{\rho^{(n_r+1)\delta}}\right) = \frac{\Lambda_0(p^*)}{\Phi'(p^*)}\int\limits_{(0,1)}\mathbb E\left(\sum_{k\in\N}|\Pi_{k}(s)|^{1+p^*}\int\limits_\mathcal P\int\limits_{(0,1]}\rho^{p^*}\underline\phi_\delta\left(\rho,\pi\right)\dd\rho\mu(\dd\pi)\right)\dd s
\]
as well as
\[
\lim_{t\to\infty}\left(\rho^{n_r\delta(1+p^*)}Z^{\bar\phi_\delta}_{\rho^{n_r\delta}}\right) = \frac{\Lambda_0(p^*)}{\Phi'(p^*)}\int\limits_{(0,1)}\mathbb E\left(\sum_{k\in\N}|\Pi_{k}(s)|^{1+p^*}\int\limits_\mathcal P\int\limits_{(0,1]}\rho^{p^*}\bar\phi_\delta\left(\rho,\pi\right)\dd\rho\mu(\dd\pi)\right)\dd s.
\]
Moreover, since  $\phi$ has $\mathbb P$-a.s. c\`adl\`ag paths, we infer that $\phi$ is continuous Lebesgue-almost everywhere. Hence, for each $\pi\in\mathcal P$ the map 
\[
\eta\mapsto\phi\left(\frac{\eta}{|\Pi_1(t_i-)|},\pi\right)
\] 
is $\mathbb P$-a.s. continuous at Lebesgue-almost every $\eta\in(0,1)$.
Resorting to the DCT, we thus infer that
\begin{align*}
& \lim_{\delta\downarrow0}\int_{(0,1)}\mathbb E\left(\sum_{k\in\N}|\Pi_{k}(s)|^{1+p^*}\int_\mathcal P\int_{(0,1]}\rho^{p^*}\underline\phi_\delta\left(\rho,\pi\right)\dd\rho\mu(\dd\pi)\right)\dd s
\\[0.5ex]
&= \int_{(0,1)}\mathbb E\left(\sum_{k\in\N}|\Pi_{k}(s)|^{1+p^*}\int_\mathcal P\int_{(0,1]}\rho^{p^*}\phi\left(\rho,\pi\right)\dd\rho\mu(\dd\pi)\right)\dd s
\\[0.5ex]
&= \lim_{\delta\downarrow0}\int_{(0,1)}\mathbb E\left(\sum_{k\in\N}|\Pi_{k}(s)|^{1+p^*}\int_\mathcal P\int_{(0,1]}\rho^{p^*}\bar\phi_\delta\left(\rho,\pi\right)\dd\rho\mu(\dd\pi)\right)\dd s.
\end{align*}
Consequently, by means of (\ref{e.t.SLLN.fprc.II.1}) 
we conclude that
\begin{align}\label{e.almost_sure_convergence}
\limsup_{\eta\downarrow0}\left(\eta^{1+p^*}Z^\phi_{\eta}\right) &= \limsup_{r\to\infty}\left(\rho^{r(1+p^*)}Z^\phi_{\rho^r}\right)\notag
\\[0.5ex]
&= \frac{\Lambda_0(p^*)}{\Phi'(p^*)}\int_{(0,1)}\mathbb E\left(\sum_{k\in\N}|\Pi_{k}(s)|^{1+p^*}\int_\mathcal P\int_{(0,1]}\rho^{p^*}\phi\left(\rho,\pi\right)\dd\rho\mu(\dd\pi)\right)\dd s
\end{align}
holds $\mathbb P$-almost surely. 

\underline{Part III} Recall that in Lemma~\ref{l.B} we showed that
\[
\lim_{\eta\downarrow0}\mathbb E\left(\eta^{1+p^*}Z^{\phi}_{\eta}\right)=\frac{1}{\Phi'(p^*)}\int_{(0,1)}\mathbb E\left(\sum_{k\in\N}|\Pi_{k}(s)|^{1+p^*}\int_\mathcal P\int_{(0,1]}\rho^{p^*}\phi\left(\rho,\pi\right)\dd\rho\mu(\dd\pi)\right)\dd s.
\]
In view of $\mathbb E(\Lambda_0(p^*))=1$ and the almost sure convergence shown in (\ref{e.almost_sure_convergence}) the corresponding convergence in $\mathscr L^1(\mathbb P)$ thus follows from Lemma~21.6 in \cite{Bau01}. 

\underline{Part IV} So far we assumed that the fragmentation process $\Pi$ is homogenous, i.e. $\alpha=0$. It remains to prove the assertion for self-similar fragmentation processes with index $\alpha\ne0$.  To this end, recall the process $(\lambda_\eta)_{\eta\in(0,1]}$ and for any $\eta\in(0,1]$ and $k\in\N$ let $\sigma_{\eta,k}$ be the {\it time of creation} of $\lambda_{\eta,k}$. That is to say, $\lambda_{\eta,k}$ corresponds to a unique block $\Pi_{m(\eta,k)}(\sigma_{\eta,k})$ and satisfies $\lambda_{\eta,k}=|\Pi_{m(\eta,k)}(\sigma_{\eta,k})|$. 
Furthermore, for any partition $\pi\in\mathcal P$ let  $\pi^+:=(\pi^+_n)_{n\in\N}$, given by   
\[
\pi^+_n:=
\begin{cases}
\pi_n, & |\pi_n|>0
\\
\emptyset, & |\pi_n|=0
\end{cases} 
\]
for every $n\in\N$, denote the subcollection of blocks with positive asymptotic frequency. Observe that 
\begin{equation}\label{e.original_def.2}
Z^\phi_\eta=\sum_{\rho\in\mathcal R}\sum_{k\in\mathcal K_\rho}\phi^{(\rho,k)}\left(\frac{\eta}{\rho},\pi(\sigma_{\rho,k})\right)
\end{equation}
holds for all $\eta\in(0,1]$, where the $\phi^{(\rho,k)}$ are independent copies of $\phi$ and where
\[
\mathcal R:=\left\{\rho\in(0,1]:\exists\,(u,m)\in\rpn\times\N\text{ such that }|\Pi_m(u)|=\rho\right\}
\]
as well as
\[
\mathcal K_\rho:=\left\{k\in\N:|\Pi_{m(\rho,k)}(\sigma_{\rho,k}-)|=\rho,\,\min\left(\bigcup_{n\in\N}\pi_n^+(\sigma_{\rho,k})\cap\Pi_{m(\rho,k)}(\sigma_{\rho,k}-)\right)\in\Pi_{m(\rho,k)}(\sigma_{\rho,k})\right\}.
\]
Note that $\Pi$ having  countably many jumps results in $\mathcal R$ being countably infinite. Moreover, notice that the representation of $Z^\phi_\eta$   in (\ref{e.original_def.2}) does not at all depend on the time parameter of $\Pi$. Indeed, all that is needed for the definition of $Z^\phi_\eta$ in (\ref{e.original_def.2}) are the sizes which the blocks of $\Pi$ attain, but not the jump times of $\Pi$. Since, according to Theorem~3~(i) of \cite{105}, any self-similar fragmentation process is a time-changed homogenous fragmentation process, the self-similar case follows from the homogenous case that we proved in Parts II and III above. Hence, the proof of Theorem~\ref{t.2.1} is complete.
\hfill$\square$

\section{Proof of Proposition~\ref{t.appendix}}\label{s.app}
The present section is devoted to the proof of Proposition~\ref{t.appendix}. Throughout this section assume that  $\Pi$ is homogenous and that $\phi$ is nonnegative. We start with some auxiliary results.

For each $t\in\rpn$ let 
\[
\lambda_1(t):=\sup_{k\in\N}|\Pi_k(t)|
\]
be the asymptotic frequency of the largest block at time $t$. In addition, for every $\eta\in(0,1]$ set
\[
\sigma_\eta:=\inf\left\{t\in\rpn:\lambda_1(t)<\eta\right\}.
\]
The following lemma provides us with an almost sure estimate of $\sigma_\eta$ by some integrable random variable.
\begin{lemma}\label{l.rmpm}
There exists some random variable $\hat\sigma\in\mathscr L^1(\mathbb P)$ such that
\[
\sigma_\eta\le \hat\sigma\cdot(-\ln(\eta)\lor1)
\]
holds $\mathbb P$-a.s. for all $\eta\in(0,1]$.
\end{lemma}

\begin{pf}
Let $\hat p\in(p^*,\bar p)$ and let $\tau:\Omega\to\rp$ be some $\mathbb P$-a.s. stopping time. Further, recall that
\[
\sup_{s\in\rpn}M_s(\hat p)\ge M_{\tau-}(\hat p)=e^{\Phi(\hat p)\tau}\sum_{k\in\N}|\Pi_k(\tau-)|^{1+\hat p}\ge e^{\Phi(\hat p)\tau}\lambda^{1+\hat p}_1(\tau-)
\]
$\mathbb P$-almost surely. Hence, observe that
\begin{equation}\label{e.rmpm.1}
\frac{-\ln(\lambda_1(\tau-))}{\tau}\ge\frac{\Phi(\hat p)}{1+\hat p}-\frac{\ln\left(\sup_{s\in\rpn}M_s(\hat p)\right)}{\tau(1+\hat p)}
\end{equation}
$\mathbb P$-almost surely. Moreover, let $\epsilon\in(0,\nicefrac{\Phi(\hat p)}{(1+\hat p)})$ and notice that
\begin{align*}
\frac{\Phi(\hat p)}{1+\hat p}-\frac{\ln\left(\sup_{s\in\rpn}M_s(\hat p)\right)}{\tau(1+\hat p)}\ge\epsilon &\Longleftrightarrow \frac{\Phi(\hat p)-\epsilon(1+\hat p)}{1+\hat p}\ge\frac{\ln\left(\sup_{s\in\rpn}M_s(\hat p)\right)}{\tau(1+\hat p)}
\\[0.5ex]
&\Longleftrightarrow \tau\ge\frac{\ln\left(\sup_{s\in\rpn}M_s(\hat p)\right)}{\Phi(\hat p)-\epsilon(1+\hat p)}
\end{align*} 
Hence, it follows from (\ref{e.rmpm.1}) that 
\[
\frac{\tau}{-\ln(\lambda_1(\tau-))}\le\frac{(1+\hat p)\tau}{\tau\Phi(\hat p)-\ln\left(\sup_{s\in\rpn}M_s(\hat p)\right)}\le\frac{1}{\epsilon}
\]
$\mathbb P$-a.s. if $\tau\ge(\ln(\sup_{s\in\rpn}M_s(\hat p)))(\Phi(\hat p)-\epsilon(1+\hat p))^{-1}$. Consequently, we infer that
\[
\frac{\sigma_\eta}{-\ln(\eta)}\le\frac{\sigma_\eta}{-\ln(\lambda_1(\sigma_\eta-))}\le\frac{(1+\hat p)\sigma_\eta}{\sigma_\eta\Phi(\hat p)-\ln\left(\sup_{s\in\rpn}M_s(\hat p)\right)}\le\frac{1}{\epsilon}.
\]
$\mathbb P$-a.s. if $\sigma_\eta\ge\ln(\sup_{s\in\rpn}M_s(\hat p))(\Phi(\hat p)-\epsilon(1+\hat p))^{-1}$. Otherwise, we have 
\[
\sigma_\eta<\frac{\ln\left(\sup_{s\in\rpn}M_s(\hat p)\right)}{\Phi(\hat p)-\epsilon(1+\hat p)}.
\]
Therefore, we obtain
\[
\sigma_\eta\le \hat\sigma\cdot(-\ln(\eta)\lor1)
\]
$\mathbb P$-a.s. for every $\eta\in(0,1]$, where 
\[
\hat\sigma:=\frac{1}{\epsilon}\lor\frac{\ln\left(\sup_{s\in\rpn}M_s(\hat p)\right)}{\Phi(\hat p)-\epsilon(1+\hat p)}.
\]
According to Proposition~3.5 in \cite{Kno11} we have 
\[
\ln\left(\sup_{s\in\rpn}M_s(\hat p)\right)\le\sup_{s\in\rpn}M_s(\hat p)\in\mathscr L^1(\mathbb P)
\] 
and thus $\hat\sigma\in\mathscr L^1(\mathbb P)$. 
\end{pf}

Recall that in Lemma~\ref{l.B} we showed that the mapping $\eta\mapsto\mathbb E(\eta^{1+p^*}Z^{\phi}_{\eta})$ has a limit as $\eta\downarrow0$. In view of the following lemma we can thus extend it to obtain a continuous function on the compact interval $[0,1]$. We shall make use of this later on.

\begin{lemma}\label{l.continuity.1}
Let $\phi$ be a random characteristic that satisfies (\ref{e.assumption_1}). Then the map $\eta\mapsto\mathbb E(\eta^{1+p^*}Z^{\phi}_{\eta})$ is continuous on $(0,1]$.
\end{lemma}

\begin{pf}
For any $k\in\N$ let $B_k(t)$ denote the block at time $t\in\rpn$ that contains the element $k$ and note that $B_1(t)=\Pi_1(t)$. By exchangeability we then infer from $(-\ln(B_1(t)))_{t\in\rpn}$ being a subordinator that also $(-\ln(B_k(t)))_{t\in\rpn}$ is a subordinator for each $k\in\N$. In addition, note that $\phi(\cdot,\pi)$, $\pi\in\mathcal P$, being c\`adl\`ag implies that $\phi(\cdot,\pi)$ has at most countably many discontinuities. 
Therefore, it follows from the proposition in \cite{Mil95} (alternatively, see \cite{HW42}, where that result was first proven) that for any $k\in\N$, $t\in\rp$ the distribution of $|B_k(t)|$ has no atoms in $[0,1]$.   Hence, for every $\eta\in(0,1)$, $\pi\in\mathcal P$, $k\in\N$ and $t\in\rpn$ the following holds:
\begin{equation}\label{e.cont_eta}
\text{The function $\phi(\cdot,\pi)$ is $\mathbb P$-a.s. continuous at } \frac{\eta}{|B_k(t)|}.
\end{equation}

Observe that the compensation formula for Poisson point processes yields that
\begin{align}\label{e.l.continuity.1.0}
& \mathbb E\left(\sum_{i\in\mathcal I}\left(\eta^{1+p^*}\phi^{(i)}\left(\frac{\eta}{|\Pi_{\kappa(t_i)}(t_i-)|},\pi(t_i)\right)-(\eta+h)^{1+p^*}\phi^{(i)}\left(\frac{\eta+h}{|\Pi_{\kappa(t_i)}(t_i-)|},\pi(t_i)\right)\right)\right)\notag
\\[0.5ex]
&\le \mathbb E\left(\int_{\rpn}\int_\mathcal P\sum_{k\in\N}\left|\eta^{1+p^*}\phi\left(\frac{\eta}{|\Pi_k(t)|},\pi\right)-(\eta+h)^{1+p^*}\phi\left(\frac{\eta+h}{|\Pi_k(t)|},\pi\right)\right|\mu(\dd\pi)\dd t\right)\notag
\\[0.5ex]
&\le  \mathbb E\left(\int_{\rpn}\int_\mathcal P\sum_{k\in\N}\eta^{1+p^*}\left|\phi\left(\frac{\eta}{|\Pi_k(t)|},\pi\right)-\phi\left(\frac{\eta+h}{|\Pi_k(t)|},\pi\right)\right|\mu(\dd\pi)\dd t\right)
\\[0.5ex]
&\qquad +\mathbb E\left(\int_{\rpn}\int_\mathcal P\sum_{k\in\N}\phi\left(\frac{\eta+h}{|\Pi_k(t)|},\pi\right)\left|\eta^{1+p^*}-(\eta+h)^{1+p^*}\right|\mu(\dd\pi)\dd t\right)\notag
\end{align}

Fix some $\eta\in(0,1)$ and $\delta_0\in(0,\eta\land(1-\eta))$. Further, let $\alpha>0$. Since $\phi(\rho,\pi)=0$ for any $\rho>1$ and $\pi\in\mathcal P$, we deduce in view of (\ref{e.assumption_1}) and Lemma~\ref{l.rmpm} that
\begin{align}\label{e.l.continuity.1.0b}
& \mathbb E\left(\int_{\rpn}\int_\mathcal P\sum_{k\in\N}\eta^{1+p^*}\sup_{h\in(-\delta,\delta)}\left|\phi\left(\frac{\eta}{|\Pi_k(t)|},\pi\right)-\phi\left(\frac{\eta+h}{|\Pi_k(t)|},\pi\right)\right|\mu(\dd\pi)\dd t\right)\notag
\\[0.5ex]
&\le 2\mathbb E\left(\int_{\left(0,\,\sigma_{\eta-\delta}\right)}\int_\mathcal P\sum_{k\in\N}\mathds1_{\{|\Pi_k(t)|\ge\eta-\delta\}}\sup_{\rho\in(\eta-\delta,1]}\left(\rho^{1+p^*}\phi\left(\rho,\pi\right)\right)\mu(\dd\pi)\dd t\right)\notag
\\[0.5ex]
&\le 2(\eta-\delta)^{-(1+\beta)}\mathbb E\left(\sigma_{\eta-\delta}\right)\mathbb E\left(\int_\mathcal P\sup_{\rho\in(\eta-\delta,1]}\left(\rho^{(1+p^*+\beta)}\phi\left(\rho,\pi\right)\right)\mu(\dd\pi)\right)
\\[0.5ex]
&< \infty\notag
\end{align}
holds for all $\delta\in(0,\delta_0)$. Thus, by means of (\ref{e.cont_eta}) and Tonelli's theorem, we can apply the DCT to infer that
\begin{align}\label{e.l.continuity.1.1}
& \lim_{h\to0}\mathbb E\left(\int_{\rpn}\int_\mathcal P\sum_{k\in\N}\eta^{1+p^*}\left|\phi\left(\frac{\eta}{|\Pi_k(t)|},\pi\right)-\phi\left(\frac{\eta+h}{|\Pi_k(t)|},\pi\right)\right|\mu(\dd\pi)\dd t\right)\notag
\\[0.5ex]
&= \int_{\rpn}\int_\mathcal P\sum_{k\in\N}\mathbb E\left(\eta^{1+p^*}\lim_{h\to0}\left|\phi\left(\frac{\eta}{|B_k(t)|},\pi\right)-\phi\left(\frac{\eta+h}{|B_k(t)|},\pi\right)\right|\mathds1_{\{k=\min(B_k(t))\}}\right)\mu(\dd\pi)\dd t\notag
\\[0.5ex]
&=0.
\end{align}

Moreover, by means of an argument as in (\ref{e.l.continuity.1.0b}) we deduce that
\begin{align}\label{e.l.continuity.1.2}
& \limsup_{h\to0}\mathbb E\left(\int_{\rpn}\int_\mathcal P\sum_{k\in\N}\phi\left(\frac{\eta+h}{|\Pi_k(t)|},\pi\right)\left|\eta^{1+p^*}-(\eta+h)^{1+p^*}\right|\mu(\dd\pi)\dd t\right)\notag
\\[0.5ex]
&\le \lim_{h\to0}\left|\eta^{1+p^*}-(\eta+h)^{1+p^*}\right|\limsup_{h\to0}\mathbb E\left(\int_\mathcal P\int_{\rpn}\sum_{k\in\N}\phi\left(\frac{\eta+h}{|\Pi_k(t)|},\pi\right)\dd t\mu(\dd\pi)\right)\notag
\\[0.5ex]
&=0.
\end{align}
Consequently, in view of (\ref{e.l.continuity.1.1}) and (\ref{e.l.continuity.1.2}), the estimate in (\ref{e.l.continuity.1.0}) results in
\begin{align*}
&  \mathbb E\left(\eta^{1+p^*}Z^{\phi}_{\eta}\right)-\lim_{h\downarrow0}\mathbb E\left((\eta+h)^{1+p^*}Z^{\phi}_{\eta+h}\right)
\\[0.5ex]
&= \lim_{h\downarrow0}\mathbb E\left(\sum_{i\in\mathcal I}\left(\eta^{1+p^*}\phi^{(i)}\left(\frac{\eta}{|\Pi_{\kappa(t_i)}(t_i-)|},\pi(t_i)\right)-(\eta+h)^{1+p^*}\phi^{(i)}\left(\frac{\eta+h}{|\Pi_{\kappa(t_i)}(t_i-)|},\pi(t_i)\right)\right)\right)
\\[0.5ex]
& =0,
\end{align*}
which completes the proof.
\end{pf}


\begin{lemma}\label{l.t.SLLN.fprc.1b.new.1b.0}
Let $\phi$ be such that (\ref{e.assumption_1}) holds. Then we have
\[
\lim_{\eta\downarrow0}\left(\eta^{1+p^*+\beta}Z^{\phi}_{\eta}\right)=0
\]
$\mathbb P$-a.s. for all $\beta>0$.
\end{lemma}

\begin{pf}

Let $a,\beta>0$ and $\rho\in(0,1)$.  Further, set
\[
Y^\phi_{a,\beta}:=\sum_{i\in\mathcal I}\mathds1_{\{t_i\le a\}}\sup_{\eta\in(0,1]}\left(\eta^{1+p^*+\beta}\phi^{(i)}\left(\eta,\pi(t_i)\right)\right),
\]
where the $\phi^{(i)}$  are independent copies of $\phi$. 
Note that
\begin{equation}\label{Y_finite}
\mathbb E\left(Y^\phi_{a,\beta}\right)<\infty.
\end{equation}
Indeed, recall that  $\phi$ and $\Pi$ are independent.  Hence, the compensation formula for Poisson point processes and Tonelli's theorem yield that
\begin{align*}
& \mathbb E\left(\sum_{i\in\mathcal I}\mathds1_{\{t_i\le a\}}\sup_{\eta\in(0,1]}\left(\eta^{1+p^*+\beta}\phi^{(i)}\left(\eta,\pi(t_i)\right)\right)\right)
\\[0.5ex]
&= a\int_\mathcal P\mathbb E\left(\sup_{\eta\in(0,1]}\left(\eta^{1+p^*+\beta}\phi\left(\eta,\pi\right)\right)\right)\mu(\dd\pi)
\\[0.5ex]
&< \infty,
\end{align*}
where the $\phi^{(i)}$ are independent copies of $\phi$, which proves that (\ref{Y_finite}) holds.

Moreover, observe that  according to Proposition~\ref{l.t.SLLN.fprc.1} and (\ref{e.l.t.SLLN.fprc.1.0}) we have $T^{(a)}_{\rho^k}<\infty$ $\mathbb P$-a.s. for ever $k\in\N_0$ and 
\[
\liminf_{k\to\infty}\frac{T^{(a)}_{\rho^k}-T^{(a)}_{\rho^{k-1}}}{T^{(a)}_{\rho^{k-1}}}\ge\liminf_{k\to\infty}\frac{T_{\rho^{k-1},\rho}}{T^{(a)}_{\rho^{k-1}}}>0
\] 
respectively. Hence, since
\[
\text{card}\left(\mathcal J^{(a)}_{\rho^k}\setminus\mathcal J^{(a)}_{\rho^{k-1}}\right)=T^{(a)}_{\rho^k}-T^{(a)}_{\rho^{k-1}},
\]
we thus infer from (\ref{Y_finite})  and a strong law of large number as in Proposition~4.1 of \cite{80} that
\begin{equation}\label{e.l.t.SLLN.fprc.1b.new.1.alpha}
\frac{\sum_{j\in\mathcal J^{(a)}_{\rho^k}\setminus\mathcal J^{(a)}_{\rho^{k-1}}}Y^{(j)}_{a,\beta}}{T^{(a)}_{\rho^k}-T^{(a)}_{\rho^{k-1}}}\to\mathbb E(Y^\phi_{a,\beta})
\end{equation}
$\mathbb P$-a.s.  as $k\to\infty$, where  the $Y^{(j)}_{a,\beta}$ are independent copies of $Y^\phi_{a,\beta}$. Furthermore, by means of the fragmentation property we have that
\begin{align*}
\rho^{l(1+p^*+2\beta)}Z^{\phi}_{\rho^l}
&= \sum_{j\in\mathcal J^{(a)}_{\rho^l}}\sum_{i\in\mathcal I}\mathds1_{\{t_{i,j}\le a\}}\rho^{l(1+p^*+2\beta)}\phi^{(i,j)}\left(\frac{\rho^l}{|\Pi_{\kappa(t_{i,j})}(t_{i,j}-)|},\pi(t_{i,j})\right)\notag
\\[0.5ex]
&= \sum_{k=1}^l\sum_{j\in\mathcal J^{(a)}_{\rho^{k}}\setminus\mathcal J^{(a)}_{\rho^{k-1}}}\sum_{i\in\mathcal I}\mathds1_{\{t_{i,j}\le a\}}\rho^{l(1+p^*+2\beta)}\phi^{(i,j)}\left(\frac{\rho^l}{|\Pi_{\kappa(t_{i,j})}(t_{i,j}-)|},\pi(t_{i,j})\right)\notag
\\[0.5ex]
&\le \rho^{-(1+p^*+\beta)}\sup_{u\in(0,1]}\left(u^{1+p^*}T^{(a)}_u\right)\rho^{l\beta}\sum_{k=1}^l\rho^{k\beta}\frac{\sum_{j\in\mathcal J^{(a)}_{\rho^{k}}\setminus\mathcal J^{(a)}_{\rho^{k-1}}}Y^{(j)}_{a,\beta}}{T^{(a)}_{\rho^{k}}-T^{(a)}_{\rho^{k-1}}}.\notag
\end{align*}
for all $l\in\N$. Consequently, resorting to (\ref{e.l.t.SLLN.fprc.1b.new.1.alpha}) as well as Proposition~\ref{l.t.SLLN.fprc.1} and bearing in mind the convergence of the geometric series $\sum_{k\in\N}\rho^{\beta k}<\infty$, we  deduce that
\[
\rho^{l(1+p^*+2\beta)}Z^{\phi}_{\rho^l}\to0
\]
$\mathbb P$-a.s. as $l\to\infty$, which completes the proof.
\end{pf}

Recall the definition of $\iota$ in (\ref{e.c.t.SLLN.1.2b.1}). By means of the extended fragmentation property we obtain that
\begin{equation}\label{e.c.t.SLLN.1.2b}
\eta^{s\iota(1+p^*)} Z^{\phi_{\iota,\eta,s}}_{\eta^{s\iota}} = \sum_{k\in\N}\lambda_{\eta,k}^{1+p^*}u_k^{1+p^*}Z^{\phi,k}_{u_k}\left|_{u_k=\frac{\eta^{s\iota}}{\lambda_{\eta,k}}}\right.
\end{equation}
holds $\mathbb P$-a.s. for any $s>1$ and $t\in\rpn$, where the $Z^{\phi,k}$ are independent copies of $Z^\phi$. Furthermore, for any $\eta\in(0,1]$ and $s>1$ set
\[
\mathcal J_{\eta,s}:=\left\{k\in\N:\lambda_{\eta,k}\ge\eta^s\right\}\quad\text{ as well as }\quad\mathcal J^\complement_{\eta,s}:=\left\{k\in\N:\lambda_{\eta,k}<\eta^s\right\}
\]
Our approach, which is inspired by the proof of Theorem~1 in \cite{HKK10}, to prove Proposition~\ref{t.appendix} is to show that asymptotically, as $\eta\downarrow0$, the conditional expectation $\mathbb E(\eta^{s\iota(1+p^*)}Z^{\phi_{\iota,\eta,s}}_{\eta^{s\iota}}|\mathscr H_\eta)$, for $s\in(1,\infty)$ sufficiently large, is a good approximation for both $n^{-\delta s\iota} Z^{\phi_{\iota,s}}_{n^{-\delta s\iota}}$ and $M_\infty(p^*)\mathbb E(\eta^{1+p^*}Z^\phi_\eta)$. For this purpose we shall need the following lemma:

\begin{lemma}\label{l.t.appendix.2}
Let $\phi$ be such that (\ref{e.assumption_1}) holds. Then  there exists some $s_0\in(1,\infty)$ such that
\[
\lim_{\eta\downarrow0}\mathbb E\left(\left.\eta^{s\iota(1+p^*)}Z^{\phi_{\iota,\eta,s}}_{\eta^{s\iota}}\right|\mathscr H_\eta\right)=\Lambda_0(p^*)\lim_{u\downarrow0}\mathbb E\left(u^{1+p^*}Z^\phi_{u}\right)
\]
$\mathbb P$-a.s. for all $s> s_0$.
\end{lemma}

\begin{pf}
For any $t\in\rpn$ and $s\in(1,\infty)$ we infer from (\ref{e.c.t.SLLN.1.2b}), in conjunction with the MCT for conditional expectations, that
\begin{align}\label{e.l.t.appendix.2.1}
& \mathbb E\left(\left.\eta^{s\iota(1+p^*)}Z^{\phi_{\iota,\eta,s}}_{\eta^{s\iota}}\right|\mathscr H_\eta\right) 
\\[0.5ex]
&= \sum_{k\in\mathcal J_{\eta,s}}\lambda^{1+p^*}_{\eta,k}\left.\mathbb E\left(u_k^{1+p^*}Z^{\phi_{\iota,\eta,s}}_{u_k}\right)\right|_{u_k=\frac{\eta^{s\iota}}{\lambda_{\eta,k}}}+\sum_{k\in\mathcal J_{\eta,s}^\complement}\lambda_{\eta,k}^{1+p^*}\left.\mathbb E\left(u_k^{1+p^*}Z^{\phi,k}_{u_k}\right)\right|_{u_k=\frac{\eta^{s\iota}}{\lambda_{\eta,k}}}\notag
\end{align}
holds $\mathbb P$-almost surely.

Since, according to  Lemma~\ref{l.B} and Lemma~\ref{l.continuity.1} we have $\sup_{u\in(0,1]}\mathbb E(u^{1+p^*}Z^{\phi}_{u})<\infty$, it follows from Lemma~3 of \cite{HKK10} that there exists some $s^*\in(1,\infty)$ such that
\begin{equation}\label{e.partI}
\lim_{\eta\downarrow0}\sum_{k\in\mathcal J_{\eta,s}^\complement}\lambda_{\eta,k}^{1+p^*}\left.\mathbb E\left(u_k^{1+p^*}Z^\phi_{u_k}\right)\right|_{u_k=\frac{\eta^{s\iota}}{\lambda_{\eta,k}}}\le\sup_{u\in(0,1]}\mathbb E\left(u^{1+p^*}Z^{\phi}_{u}\right)\lim_{\eta\downarrow0}\sum_{k\in\mathcal J_{\eta,s}^\complement}\lambda_{\eta,k}^{1+p^*}=0
\end{equation}
$\mathbb P$-a.s. for all $s> s^*$. 
Let us now consider the first summand on the right-hand side of (\ref{e.l.t.appendix.2.1}). We aim at showing that there exists some $s^{**}\in(1,\infty)$ such that
\begin{equation}\label{e.l.t.appendix.2.1b}
\sum_{k\in\mathcal J_{\eta,s}}\lambda^{1+p^*}_{\eta,k}\left.\mathbb E\left(u_k^{1+p^*}Z^{\phi_{\iota,\eta,s}}_{u_k}\right)\right|_{u_k=\frac{\eta^{s\iota}}{\lambda_{\eta,k}}}\to \Lambda_0(p^*)\lim_{u\downarrow0}\mathbb E\left(u^{1+p^*}Z^\phi_{u}\right)
\end{equation} 
holds $\mathbb P$-a.s. for every $s\ge s^{**}$ as $\eta\downarrow0$. 

For the time being, fix some $s\in(1,\infty)$. 
Further, note that Lemma~\ref{l.B} together with Lemma~\ref{l.continuity.1} yields that   the maps $g_0:[0,1]\to\rpn$ and $g_\eta:[0,1]\to\rpn$, given by
\[
g_0(0)=\lim_{\rho\downarrow0}\mathbb E\left(\rho^{1+p^*}Z^\phi_{\rho}\right)\qquad\text{and}\qquad g_\eta=\lim_{\rho\downarrow0}\mathbb E\left(\rho^{1+p^*}Z^{\phi_{\iota,\eta,s}}_{\rho}\right)
\]
as well as
\[
g_0(u)=\mathbb E\left(u^{1+p^*}Z^\phi_{u}\right)\qquad\text{and}\qquad g_\eta(u)=\mathbb E\left(u^{1+p^*}Z^{\phi_{\iota,\eta,s}}_{u}\right),
\]
are continuous  for all $\eta\in(0,1]$ and $u\in(0,1]$. Moreover, since 
\[
\forall\,\eta,u\in(0,1],\,\forall\,\gamma\le\eta u^{-s(\iota-1)}: \phi_{\iota,\eta,s}\left(\frac{u}{\gamma}\right)=\phi\left(\frac{u}{\gamma}\right),
\] 
we deduce that  
\[
\forall\,u\in(0,1],\,\forall\,\eta\le u^{\frac{1}{s(\iota-1)}},\,\forall\pi\in\mathcal P,\,\forall\,k\in\N:\phi_{\iota,\eta,s}\left(\frac{u}{|\pi_k|}\right)=\phi\left(\frac{u}{|\pi_k|}\right).
\] 
Therefore, for every $u\in[0,1]$ we have that $g_\eta(u)\to g_0(u)$ as $\eta\downarrow0$. Consequently, since for each $u\in[0,1]$ the mapping $\eta\mapsto g_\eta(u)$ is nonincreasing, we obtain by resorting to Dini's theorem that
\[
g_\eta(u)\to g_0(u) 
\]
uniformly in $u\in(0,1]$ as $\eta\downarrow0$, i.e.
\[
\sup_{u\in(0,1]}\left|u^{1+p^*}\mathbb E\left(Z^{\phi_{\iota,\eta,s}}_{u}-Z^{\phi}_{u}\right)\right|\to0
\]
as $\eta\downarrow0$. Hence,
\begin{align}\label{e.unif_conv}
& \left|\lim_{\eta\downarrow0}\left.\mathbb E\left(u_k^{1+p^*}Z^{\phi_{\iota,\eta,s}}_{u_k}\right)\right|_{u_k=\frac{\eta^{s\iota}}{\lambda_{\eta,k}}}-\lim_{\eta\downarrow0}\left.\mathbb E\left(u_k^{1+p^*}Z^\phi_{u_k}\right)\right|_{u_k=\frac{\eta^{s\iota}}{\lambda_{\eta,k}}}\right| \notag
\\[0.5ex]
&\le \lim_{\eta\downarrow0}\sup_{u\in\rpn}\left|u^{1+p^*}\mathbb E\left(Z^{\phi_{\iota,\eta,s}}_{u}-Z^\phi_{u}\right)\right|
\\[0.5ex]
& =0\notag
\end{align}
$\mathbb P$-a.s. for all $\eta\in(0,1]$, $k\in\mathcal J_{\eta,s}$ and $s\in(1,\infty)$.
Furthermore, since
\[
\frac{\eta^{s\iota}}{\lambda_{\eta,k}}\le\eta^{s(\iota-1)} 
\]
holds for all $\eta\in(0,1]$, $k\in\mathcal J_{\eta,s}$ and $s\in(1,\infty)$, we have
\[
\left|\left.\mathbb E\left(u_k^{1+p^*}Z^\phi_{u_k}\right)\right|_{u_k=\frac{\eta^{s\iota}}{\lambda_{\eta,k}}}
-\lim_{u\downarrow0}\mathbb E\left(u^{1+p^*}Z^\phi_{u}\right)\right|\to0 
\]
uniformly in $k\in\mathcal J_{\eta,s}$ as $\eta\downarrow0$. Thus we deduce from (\ref{e.unif_conv}) that
\begin{equation}\label{e.l.t.appendix.2.1e}
\left.\mathbb E\left(u_k^{1+p^*}Z^{\phi_{\iota,\eta,s}}_{u_k}\right)\right|_{u_k=\frac{\eta^{s\iota}}{\lambda_{\eta,k}}}-\lim_{u\downarrow0}\mathbb E\left(u^{1+p^*}Z^\phi_{u}\right)\to0 
\end{equation}
uniformly in $k\in\mathcal J_{\eta,s}$ as $\eta\downarrow0$.
Moreover, according to Lemma~3 of \cite{HKK10}  there exists some $s^{**}\in(1,\infty)$ such that 
\[
\lim_{\eta\downarrow0}\sum_{k\in\mathcal J_{\eta,s}}\lambda^{1+p^*}_{\eta,k}=\Lambda_0(p^*)
\]
holds $\mathbb P$-a.s. for all $s\ge s^{**}$. In conjunction with (\ref{e.l.t.appendix.2.1e}) this proves (\ref{e.l.t.appendix.2.1b}).

Consequently, choosing $s_0:=s^*\lor s^{**}$ the assertion of the lemma follows from (\ref{e.partI}) and (\ref{e.l.t.appendix.2.1b}).
\end{pf}

Now we are ready to prove Proposition~\ref{t.appendix}.

{\bf Proof of Proposition~\ref{t.appendix}}
The proof is divided into three parts. The first two parts are concerned with establishing an $\mathscr L^2$-estimate of the difference between the counted process and its conditional expectation. In the third part we use this estimate to prove the desired almost sure convergence.

\underline{Part I} As in Lemma~1 of \cite{83} an application of Fatou's lemma (for conditional expectations) results in
\begin{equation}\label{e.l.t.appendix.3.1}
\mathbb E\left(\left.\left|\sum_{n\in\N}X_n\right|^2\,\right|\mathscr H_\eta\right)\le4\sum_{n\in\N}\mathbb E\left(\left.\left|X_n\right|^2\,\right|\mathscr H_t\right)
\end{equation}
for any  sequence $(X_n)_{n\in\N}$ of independent centred random variables. Moreover, according to Jensen's inequality we have
\begin{equation}\label{e.l.t.appendix.3.2}
|u+v|^2\le2\left(|u|^2+|v|^2\right)
\end{equation}
for all $u,v\in\R$.

For the time being, let $\eta\in(0,1]$ as well as $s>1$. By means of (\ref{e.c.t.SLLN.1.2b}) we obtain
\begin{align}\label{e.l.t.appendix.3.3}
& \eta^{s\iota(1+p^*)}Z^{\phi_{\iota,\eta,s}}_{\eta^{s\iota}}-\mathbb E\left(\left.\eta^{s\iota(1+p^*)}Z^{\phi_{\iota,\eta,s}}_{\eta^{s\iota}}\right|\mathscr H_\eta\right) \notag
\\[0.5ex]
&= \sum_{k\in\mathcal J_{\eta,s}}\lambda^{1+p^*}_{\eta,k}\left(Z^{(k)}-\mathbb E\left(\left.Z^{(n)}\right|\mathscr H_\eta\right)\right)\notag
\\[0.5ex]
&\qquad +\sum_{k\in\mathcal J_{\eta,s}^\complement}\lambda_{\eta,k}^{1+p^*}u_k^{1+p^*}Z^{\phi,k}_{u_k}\left|_{u_k=\frac{\eta^{s\iota}}{\lambda_{\eta,k}}}\right.-\mathbb E\left(\left.\sum_{k\in\mathcal J_{\eta,s}^\complement}\lambda_{\eta,k}^{1+p^*}u_k^{1+p^*}Z^{\phi,k}_{u_k}\left|_{u_k=\frac{\eta^{s\iota}}{\lambda_{\eta,k}}}\right.\right|\mathscr H_\eta\right)\notag
\\[0.5ex]
&= \sum_{k\in\mathcal J_{\eta,s}}\lambda^{1+p^*}_{\eta,k}\left(Z^{(k)}-\mathbb E\left(\left.Z^{(k)}\right|\mathscr H_\eta\right)\right),
\end{align}
where conditional on $\mathscr H_t$ the $Z^{(k)}$ are independent and satisfy
\begin{equation}\label{e.Zk}
\mathbb P\left(\left.Z^{(k)}\in\cdot\right|\mathscr H_\eta\right)=\left.\mathbb P\left(u_k^{1+p^*}Z^{\phi_{\iota,\eta,s}}_{u_k}\in\cdot\right)\right|_{u_k=\frac{\eta^{s\iota}}{\lambda_{\eta,k}}}
\end{equation}
$\mathbb P$-almost surely. Since $\mathbb E(u^{1+p^*}Z^{\phi}_u-\mathbb E(u^{1+p^*}Z^{\phi}_u|\mathscr H_\eta))=0$ for all $u\in(0,1]$, we can apply (\ref{e.l.t.appendix.3.1}) in order to deduce from (\ref{e.l.t.appendix.3.3}) that 
\begin{align}\label{e.l.t.appendix.3.4}
& \mathbb E\left(\left.\left|\eta^{s\iota(1+p^*)}Z^{\phi_{\iota,\eta,s}}_{\eta^{s\iota}}-\mathbb E\left(\left.\eta^{s\iota(1+p^*)}Z^{\phi_{\iota,\eta,s}}_{\eta^{s\iota}}\right|\mathscr H_\eta\right)\right|^2\,\right|\mathscr H_\eta\right)\notag
\\[0,5ex]
&\le 2^2\sum_{k\in\mathcal J_{\eta,s}}\lambda_{\eta,k}^{p(1+p^*)}\mathbb E\left(\left.\left|Z^{(k)}-\mathbb E\left(\left.Z^{(k)}\right|\mathscr H_\eta\right)\right|^2\,\right|\mathscr H_\eta\right)\notag
\\[0.5ex]
&\le 2^3\sum_{k\in\mathcal J_{\eta,s}}\lambda_{\eta,k}^{p(1+p^*)}\mathbb E\left(\left.\left(Z^{(k)}\right)^2+\mathbb E\left(\left.Z^{(k)}\right|\mathscr H_\eta\right)^2\,\right|\mathscr H_\eta\right)
\\[0.5ex]
&\le 2^3\sum_{k\in\mathcal J_{\eta,s}}\lambda_{\eta,k}^{p(1+p^*)}\mathbb E\left(\left.\left(Z^{(k)}\right)^2+\mathbb E\left(\left.\left(Z^{(k)}\right)^2\right|\mathscr H_\eta\right)\,\right|\mathscr H_\eta\right)\notag
\\[0.5ex]
&= 2^4\sum_{k\in\mathcal J_{\eta,s}}\lambda_{\eta,k}^{p(1+p^*)}\mathbb E\left(\left.\left(Z^{(k)}\right)^2\right|\mathscr H_\eta\right),\notag
\end{align}
where the $Z^{(k)}$ are the same random variables given by (\ref{e.Zk}) that appear in (\ref{e.l.t.appendix.3.3}). Notice that the first estimate in (\ref{e.l.t.appendix.3.4}) results from (\ref{e.l.t.appendix.3.1}) and (\ref{e.l.t.appendix.3.3}), and the second estimate holds by means of (\ref{e.l.t.appendix.3.2}). The third estimate is a consequence of Jensen's inequality for conditional expectations.

\underline{Part II} For the remainder of the proof  fix some $s\in(1,\infty)$ as well as $\epsilon\in(0,1)$ and choose some $\beta>0$ such that $\beta<(1+p^*)(s\iota-1)^{-1}$. Furthermore,  for every $\eta,\rho\in(0,1]$ define
\[
\tilde Z^{\phi_{\iota,\eta,s}}_\rho:=Z^{\phi_{\iota,\eta,s}}_\rho\mathds1_{\{\rho\le\hat\rho\}},
\]
where in view of Lemma~\ref{l.t.SLLN.fprc.1b.new.1b.0} the random variable $\hat\rho:\Omega\to(0,1)$ is defined such that 
\[
\rho^{1+p^*+\beta}Z^{\phi_{\iota,\eta,s}}_{\rho}\le\epsilon
\]
holds $\mathbb P$-a.s. for all $\rho\le\hat\rho$. Hence, we have
\[
\sup_{\rho\in(0,1]}\mathbb E\left(\left(\rho^\beta\rho^{1+p^*}\tilde Z^{\phi_{\iota,\eta,s}}_{\rho}\right)^2\right)\le\epsilon^2<\epsilon.
\]
Observing that (\ref{e.l.t.appendix.3.4})  also holds with $\tilde Z^{\phi_{\iota,\eta,s}}$ instead of $Z^{\phi_{\iota,\eta,s}}$, and taking expectations on both sides of the resulting estimate, we thus obtain
\begin{align}\label{e.l.t.appendix.3.5}
& \mathbb E\left(\left|\eta^{s\iota(1+p^*)}\tilde Z^{\phi_{\iota,\eta,s}}_{\eta^{s\iota}}-\mathbb E\left(\left.\eta^{s\iota(1+p^*)}\tilde Z^{\phi_{\iota,\eta,s}}_{\eta^{s\iota}}\right|\mathscr H_\eta\right)\right|^2\right) \notag
\\[0.5ex]
&\le 32\mathbb E\left(\sum_{k\in\mathcal J_{\eta,s}}\lambda_{\eta,k}^{p(1+p^*)}\left.\mathbb E\left(\left(u_k^{1+p^*}\tilde Z^{\phi_{\iota,\eta,s}}_{u_k}\right)^2\right)\right|_{u_k=\frac{\eta^{s\iota}}{\lambda_{\eta,k}}}\right)
\\[0.5ex]
&\le 32\eta^{-\beta(s\iota-1)}\epsilon\mathbb E\left(\sum_{k\in\mathcal J_{\eta,s}}\lambda_{\eta,k}^{p(1+p^*)}\right),\notag
\end{align}
where the second inequality follows from $\eta^{s\iota}\lambda_{\eta,k}^{-1}\ge\eta^{s\iota-1}$ for all $k\in\mathcal J_{\eta,s}$ and from the obvious fact that $\tilde Z^{\phi_{\iota,\eta,s}}_{u_k}\le \tilde Z^{\phi}_{u_k}$ for each $k\in\N$.  
Moreover, note that 
\[
\mathbb E\left(\sum_{k\in\mathcal J_{\eta,s}}\lambda_{\eta,k}^{2(1+p^*)}\right)\le \eta^{1+p^*}\mathbb E\left(\Lambda_\eta(p^*)\right)=\eta^{1+p^*},
\]
since $\Lambda(p^*)$ is a unit-mean martingale. Consequently, we infer from (\ref{e.l.t.appendix.3.5})  that
\begin{equation}\label{e.l.t.appendix.3.6}
\mathbb E\left(\left|\eta^{s\iota(1+p^*)}\tilde Z^{\phi_{\iota,\eta,s}}_{\eta^{s\iota}}-\mathbb E\left(\left.\eta^{s\iota(1+p^*)}\tilde Z^{\phi_{\iota,\eta,s}}_{\eta^{s\iota}}\right|\mathscr H_\eta\right)\right|^2\right)\le32\epsilon\eta^{1+p^*-\beta(s\iota-1)}.
\end{equation}

\underline{Part III}
Recall that $\beta,\epsilon$ and $s$ are fixed as in Part II of this proof.  By means  of   the Chebyshev-Markov inequality the estimate in (\ref{e.l.t.appendix.3.6}) results in
\begin{align*}
& \sum_{k\in\N}\mathbb P\left(\left|\eta^{k\delta s\iota(1+p^*)}\tilde Z^{\phi_{\iota,\eta^{k\delta},s}}_{\eta^{k\delta s\iota}}-\mathbb E\left(\left.\eta^{k\delta s\iota(1+p^*)}\tilde Z^{\phi_{\iota,\eta^{k\delta},s}}_{\eta^{k\delta s\iota}}\right|\mathscr H_{\eta^{k\delta}}\right)\right|>\epsilon\right)
\\[0.5ex]
&\le \frac{1}{\varepsilon^2}\sum_{k\in\N}\left\|\eta^{k\delta s\iota(1+p^*)}\tilde Z^{\phi_{\iota,\eta^{k\delta},s}}_{\eta^{k\delta s\iota}}-\mathbb E\left(\left.\eta^{s\iota(1+p^*)}\tilde Z^{\phi_{\iota,\eta^{k\delta},s}}_{\eta^{k\delta s\iota}}\right|\mathscr H_{\eta^{k\delta}}\right)\right\|^2_{\mathscr L^2(\mathbb P)}
\\[0.5ex]
&\le \frac{32\epsilon}{\varepsilon^2}\sum_{k\in\N}\eta^{k\delta(1+p^*-\beta(s\iota-1))}
\\[0.5ex]
&< \infty
\end{align*}
for every $\eta\in(0,1)$ and all $\delta>0$. Hence, we infer from the Borel-Cantelli lemma that
\begin{equation}\label{e.t.appendix.1}
\lim_{k\to\infty}\left|\eta^{k\delta s\iota(1+p^*)}\tilde Z^{\phi_{\iota,\eta^{k\delta},s}}_{\eta^{k\delta s\iota}}-\mathbb E\left(\left.\eta^{k\delta s\iota(1+p^*)}\tilde Z^{\phi_{\iota,\eta^{k\delta},s}}_{\eta^{k\delta s\iota}}\right|\mathscr H_{\eta^{k\delta}}\right)\right|=0
\end{equation}
holds $\mathbb P$-a.s. for all $\eta\in(0,1)$, $\delta>0$.  Recall that $s\in(1,\infty)$ was chosen arbitrarily. In particular, we can assume without loss of generality that $s\ge s_0$, where $s_0$ is given by Lemma~\ref{l.t.appendix.2}. In view of the triangle inequality we thus deduce from Lemma~\ref{l.t.appendix.2} and (\ref{e.t.appendix.1})  that
\begin{align}\label{e.t.appendix.2}
& \lim_{k\to\infty}\left|\eta^{k\delta s\iota(1+p^*)}Z^{\phi_{\iota,{\eta}^{k\delta},s}}_{\eta^{k\delta s\iota}}-\Lambda_0(p^*)\lim_{u\downarrow0}\mathbb E\left(u^{1+p^*}Z^{\phi_{\iota,u}}_u\right)\right|\notag
\\[0.5ex]
&\le \lim_{k\to\infty}\left|\eta^{k\delta s\iota(1+p^*)}Z^{\phi_{\iota,{\eta}^{k\delta},s}}_{\eta^{k\delta s\iota}}-\eta^{k\delta s\iota(1+p^*)}\tilde Z^{\phi_{\iota,{\eta}^{k\delta},s}}_{\eta^{k\delta s\iota}}\right|\notag
\\[0.5ex]
&\qquad +\lim_{k\to\infty}\left|\eta^{k\delta s\iota(1+p^*)}\tilde Z^{\phi_{\iota,{\eta}^{k\delta},s}}_{\eta^{k\delta s\iota}}-\mathbb E\left(\left.\eta^{k\delta s\iota(1+p^*)}\tilde Z^{\phi_{\iota,{\eta}^{k\delta},s}}_{\eta^{k\delta s\iota}}\right|\mathscr H_{{\eta}^{k\delta}}\right)\right| \notag
\\[0.5ex]
&\qquad +\lim_{k\to\infty}\left|\mathbb E\left(\left.\eta^{k\delta s\iota(1+p^*)}\tilde Z^{\phi_{\iota,\eta^{k\delta},s},\epsilon}_{\eta^{k\delta s\iota}}-\eta^{k\delta s\iota(1+p^*)}Z^{\phi_{\iota,\eta^{k\delta},s}}_{\eta^{k\delta s\iota}}\right|\mathscr H_{\eta^{k\delta}}\right)\right|
\\[0.5ex]
&\qquad +\lim_{k\to\infty}\left|\mathbb E\left(\left.\eta^{k\delta s\iota(1+p^*)}Z^{\phi_{\iota,{\eta}^{k\delta},s}}_{\eta^{k\delta s\iota}}\right|\mathscr H_{{\eta}^{k\delta}}\right)-\Lambda_0(p^*)\lim_{u\downarrow0}\mathbb E\left(u^{1+p^*}Z^{\phi_{\iota,u}}_u\right)\right|\notag
\\[0.5ex]
&= 0\notag
\end{align}
$\mathbb P$-a.s. for all $\eta\in(0,1)$ and $\delta>0$. Let $\rho\in(0,1)$. Setting $\eta:=\rho^{\nicefrac{1}{s}}$ in (\ref{e.t.appendix.2}) we obtain that
\[
\lim_{k\to\infty}\rho^{k\delta\iota(1+p^*)} Z^{\phi_{\iota,\rho^{k\delta}}}_{\rho^{k\delta\iota}}= \Lambda_0(p^*)\lim_{u\downarrow0}\mathbb E\left(u^{1+p^*}Z^{\phi_{\iota,u}}_u\right)
\]
holds $\mathbb P$-a.s. for all $\delta>0$, which completes the proof of Proposition~\ref{t.appendix}.
\hfill$\square$

\end{document}